\documentclass[12pt,leqno]{article}

\usepackage{amsfonts}
\usepackage{amsmath}
\usepackage{amssymb}
\usepackage{epsf}

\newcommand{\Xcomment}[1]{}

\makeatletter
\renewcommand{\section}{\@startsection{section}{1}{0pt}%
{-3.5ex plus -1ex minus -.2ex}{2.3ex plus .2ex}%
{\normalfont\Large}}

\makeatletter \@addtoreset{equation}{section}
\makeatother

\newtheorem{theorem}{Theorem}[section]

\newtheorem{corollary}[theorem]{Corollary}
\newtheorem{prop}[theorem]{Proposition}

\newenvironment{proof}{\noindent{\bf Proof~}}%
{$\qed$\vskip 0.1in}

\newenvironment{numitem}{\refstepcounter{equation}\begin{enumerate}%
\item[(\arabic{section}.\arabic{equation})]$\quad$}{\end{enumerate}}

\newcommand{\refeq}[1]{(\ref{eq:#1})}  

\def\qed{ \ \vrule width.2cm height.2cm depth0cm}

\def\tilde{\widetilde}
\def\hat{\widehat}
\def\bar{\overline}
\def\eps{\epsilon}

\def\Rset{{\mathbb R}}
\def\Zset{{\mathbb Z}}

\def\Ascr{\mathcal{A}}

\def\Cscr{\mathcal{C}}
\def\Escr{\mathcal{E}}
\def\Fscr{\mathcal{F}}
\def\Gscr{\mathcal{G}}

\def\Iscr{\mathcal{I}}

\def\Kscr{\mathcal{K}}

\def\Mscr{\mathcal{M}}

\def\Vscr{\mathcal{V}}

\def\bfa{{\bf a}}
\def\bfb{{\bf b}}
\def\bfc{{\bf c}}
\def\bfd{{\bf d}}
\def\bfq{{\bf q}}

\def\bfzero{{\bf 0}}

\def\Pin{P^{\rm in}}
\def\Pout{P^{\rm out}}
\def\ellin{t}
\def\ellout{h}

\def\eNW{e^{\rm NW}}
\def\eSW{e^{\rm SW}}
\def\eNE{e^{\rm NE}}
\def\eSE{e^{\rm SE}}
\def\PNE{P^{\rm NE}}
\def\PSE{P^{\rm SE}}

\def\lvert{{left}}
\def\rvert{{right}}
\def\botvert{{bottom}}
\def\topvert{{top}}

\def\Int{{\rm Int}}

\def\sprod{\;\hbox{\unitlength=1mm\begin{picture}(4,4)%
\put(0,0){$\triangleright$}\put(1.5,0){$\triangleleft$}%
\end{picture}}\;}

\begin{document}

\begin{center}
{\large\bf  The crossing model for regular $A_n$-crystals}%
\footnote[1]{This research was supported in part by NWO--RFBR grant
047.011.2004.017 and by RFBR grant 05-01-02805
CNRSL\_a.}
\end{center}

\begin{center}
{\sc Vladimir~I.~Danilov}\footnote[2]
{Central Institute of Economics and Mathematics
of the RAS, 47, Nakhimovskii Prospect, 117418 Moscow, Russia;
emails: danilov@cemi.rssi.ru, koshevoy@cemi.rssi.ru.},
{\sc Alexander~V.~Karzanov}\footnote[3]{Institute for System Analysis of
the RAS, 9, Prospect 60 Let Oktyabrya, 117312 Moscow, Russia; email:
sasha@cs.isa.ru. A part of the original version of this paper was written
while this author was visiting group PNA1 at CWI, Amsterdam in fall of
2006 and was
supported by a grant from this Center.},\\
{\sc and Gleb~A.~Koshevoy}$^2$
\end{center}

\begin{quote}
{\bf Abstract.} A regular $A_n$-crystal is an edge-colored directed graph,
with $n$ colors, related to an irreducible highest weight integrable
module over $U_q(sl_{n+1})$. Based on Stembridge's local axioms for
regular simply-laced crystals and a structural characterization of regular
$A_2$-crystals in~\cite{DKK-07}, we present a new combinatorial
construction, the so-called {\em crossing model}, and prove that this
model generates precisely the set of regular $A_n$-crystals.

Using the model, we obtain a series of results on the combinatorial
structure of such crystals and properties of their subcrystals.
  \end{quote}



\baselineskip=16pt

\section{Introduction} \label{sec:intr}

The notion of a {\em crystal} introduced by Kashiwara~\cite{kas-90,kas-91} has
proved its importance in representation theory. This is an edge-colored
directed graph, with $n$ colors, in which each connected monochromatic subgraph
is a finite path, and there are certain interrelations on the lengths of such
paths, described via coefficients of an $n\times n$ Cartan matrix $M$ (this
matrix characterizes the {\em type} of a crystal). The central role in the
theory of Kashiwara is played by crystals of representations, or {\em regular}
crystals; these are associated to irreducible highest weight integrable modules
(representations) over the quantum enveloping algebra related to $M$. There are
several global models to characterize the regular crystals for a variety of
types; e.g., via generalized Young tableaux~\cite{KN-94}, Lusztig's canonical
bases~\cite{Lusztig}, Littelmann's path model~\cite{lit-95,lit-98}.

Stembridge~\cite{ste-03} pointed out a list of ``local'' graph-theoretic
defining axioms for the regular {\em simply-laced} crystals. These concern
simply-laced Cartan matrices $M$, i.e., those having coefficients
$m_{ii}=2$ and $m_{ij}=m_{ji}\in\{0,-1\}$ for $i\ne j$. He showed that if
$M$ has full rank, then for each $n$-tuple $\bfc=(c_1,\ldots,c_n)$ of
nonnegative integers, there is precisely one graph $K$ satisfying his
axioms and such that: $K$ is acyclic and has a unique minimal vertex ({\em
source}) $s$, and the lengths of maximal monochromatic paths with colors
$1,\ldots,n$ beginning at $s$ are equal to $c_1,\ldots,c_n$, respectively.
Moreover, $K$ is a regular crystal related to $M$ (it is the crystal graph
of the integrable module of highest weight $\sum_ic_i\omega_i$ over the
corresponding quantum enveloping algebra, where $\omega_i$ is $i$-th
fundamental weight). So in this case (and when $M$ is fixed) $\bfc$ may be
regarded as the {\em parameter} of $K$, and we may denote $K$ by
$K(\bfc)$.

This paper is devoted to a combinatorial study of regular simply-laced
crystals of $A_n$-type, or {\em regular $A_n$-crystals}; for brevity we
throughout call them {\em RAN-crystals}. They are related to the algebra
$U_q(sl_{n+1})$, and the off-diagonal coefficients $m_{ij}$ of the Cartan
matrix (which is of full rank) are equal to $-1$ if $|i-j|=1$, and 0
otherwise.

In our previous paper~\cite{DKK-07} we described the combinatorial
structure of regular $A_2$-crystals $K$ and demonstrated additional
combinatorial and polyhedral properties of these crystals and their
extensions. The structure turns out to be rather transparent: $K$ always
has a (unique) source, and therefore, $K=K(c_1,c_2)$ for some
$c_1,c_2\in\Zset_+$, and it can be produced by a certain operation of
replicating and gluing together from the crystals $K(c_1,0)$ and
$K(0,c_2)$. The latter crystals are of simple form and are viewed as
triangle-shaped parts of square grids. (In fact, $K$ is the largest
component of the tensor product of $K(c_1,0)$ and $K(0,c_2)$.)

When $n>2$, the structure of an RAN-crystal becomes much more sophisticated,
even for $n=3$. To explore this structure, in this paper we introduce a certain
combinatorial construction, called the {\em crossing model}. This model
consists of three ingredients: (i) a finite directed graph $G$, called the {\em
supporting graph}, depending only on the number $n$ of colors; (ii) a set
$\Fscr$ of integer-valued {\em feasible functions} on the vertices of $G$,
depending on a parameter $\bfc\in\Zset_+^n$; and (iii) $n$ sets
$\Escr_1,\ldots,\Escr_n$, each consisting of {\em transformations} $f\mapsto
f'$ of feasible functions. (In fact, the crossing model is a sort of
decomposition of the Gelfand-Tsetlin pattern model~\cite{GT-50}.)

Our main working theorem asserts that the $n$-colored directed graph formed by
$\Fscr$ as the vertex set and by $\Escr_1,\ldots,\Escr_n$ as the edge sets of
colors $1,\ldots,n$, respectively, is isomorphic to the RAN-crystal $K(\bfc)$.
In addition, we explain that {\em any} finite graph satisfying Stembridge's
axioms for the $A_n$ case has a source. Therefore, the crossing model produces
precisely the set of crystals of representations for $U_q(sl_{n+1})$. Our
construction and proofs rely merely on Stembridge's axiomatics and
combinatorial arguments and do not appeal explicitly to powerful tools, such as
the Path Model, or so.

Then we take advantages from the description of RAN-crystals via the
crossing model. The supporting graph $G$ consists of $n$ pairwise disjoint
subgraphs $G^1,\ldots,G^n$, and given a parameter $\bfc$, the values of
any feasible function to each $G^i$ ($i=1,\ldots,n$) are between 0 and
$c_i$. The feasible functions that are constant within each $G^i$ are of
especial interest to us. We refer to the vertices of the crystal $K(\bfc)$
corresponding to such functions as {\em principal} ones, and to the set
$\Pi$ of these vertices as the {\em principal lattice}. So there are
$(c_1+1)\times\ldots\times(c_n+1)$ principal vertices, each corresponding
to an $n$-tuple $\bfa=(a_1,\ldots,a_n)\in\Zset_+^n$ with $\bfa\le \bfc$,
being denoted as $v[\bfa]$. The principal lattice $\Pi$ is proved to have
the following properties:

(i) for any $\bfa,\bfb\in\Zset_+^n$ with $\bfa\le\bfb\le\bfc$, the
interval of $K(\bfc)$ between the principal vertices $v[\bfa]$ and
$v[\bfb]$ is isomorphic to the RAN-crystal $K(\bfb-\bfa)$;

(ii) there are exactly $|\Pi|$ maximal (connected) subcrystal of $K(\bfc)$
with colors $1,\ldots,n-1$ and each of them contains exactly one principal
vertex; a similar property takes place for the maximal subcrystals with
colors $2,\ldots,n$.

We also establish other features of maximal subcrystals $K'$ with colors
$1,\ldots,n-1$ (or $2,\ldots,n$). In particular, the parameter of $K'$ is
expressed by an explicit linear function of $\bfc$ and $\bfa$, where
$v[\bfa]$ is the principal vertex in $K'$.

The crossing model enables us to reveal one more interesting object in the
crystal $K(\bfc)$. When a feasible function varies within some subgraph $G^i$
and is constant within each of the other subgraphs $G^j$ of the supporting
graph $G$, we obtain an $n$-colored subcrystal of $K(\bfc)$ having the
parameter $\bfc'$ such that $c'_i=c_i$ and $c'_j=0$ for $j\ne i$. (This is the
crystal graph of the integrable module of $U_q(sl_{n+1})$ with the highest
weight $c_i\omega_i$.) The union of these subcrystals (for all $i$) forms a
canonical subgraph that we call the {\em skeleton} of $K(\bfc)$. It coincides
with the whole crystal $K(\bfc)$ when $n=2$, and is typically smaller when
$n\ge 3$.

By use of the crossing model, we also can derive natural infinite analogs
of RAN-crystals, in which some or all maximal monochromatic paths are
infinite (this generalizes the construction of infinite $A_2$-crystals
in~\cite{DKK-07}).

This paper is organized as follows. Section~\ref{sec:An} states
Stembridge's axioms for RAN-crystals, recalls some basic properties of
crystals, and briefly reviews results on $A_2$-crystals
from~\cite{DKK-07}. Also, relying on a structural characterization of
regular $A_2$-crystals, we explain in this section that any RAN-crystal
has a source (Corollary~\ref{cor:RANsource}). The crossing model is
described throughout Section~\ref{sec:model} (concerning the supporting
graph and feasible functions) and Section~\ref{sec:moves} (concerning
transformations of feasible functions). The equivalence between the
objects generated by the crossing model and the RAN-crystals is proved in
Section~\ref{sec:mod-cryst} (Theorem~\ref{tm:equiv}).
Section~\ref{sec:princ} introduces the principal lattice, principal
intervals and the skeleton of an RAN-crystal and explains relations
between these objects. Also infinite analogs of RAN-crystals and their
properties are discussed in this section. Section~\ref{sec:subcryst} is
devoted to a study of maximal $(n-1)$-colored subcrystals; here we prove
the above-mentioned relation between these subcrystals and the principal
lattice, compute their parameters and multiplicities, and discusses
additional issues.

 \medskip
Our study of RAN-crystals by use of the crossing model will be continued
in the forthcoming paper~\cite{DKK-08b} where we characterize the pairwise
intersections of maximal subcrystals with colors $1,\ldots,n-1$ and colors
$2,\ldots,n$ and, as a result, obtain a recursive description of the
combinatorial structure and an algorithm of assembling of the RAN-crystal
for a given parameter. (Also, using results on $B_2$-crystals
from~\cite{DKK-08a}, we discuss there a relation between RAN-crystals and
regular crystals of types B and C.)

\section{Axioms of RAN-crystals and backgrounds}
\label{sec:An}

Throughout, by an $n$-colored digraph we mean a (finite or infinite) directed
graph $K=(V(K),E(K))$ with vertex set $V(K)$ and with edge set $E(K)$
partitioned into $n$ subsets $E_1,\ldots,E_n$.
We say that an edge in $E_i$ has {\em color} $i$ and for brevity
call it an $i$-{\em edge}.

 \bigskip
 \noindent
 {\bf\large \arabic{section}.1. Axioms.}
Stembridge~\cite{ste-03} pointed out local graph-theoretic axioms that
precisely characterize the set of regular simply-laced crystals. The
RAN-crystals (which form a subclass of regular simply-laced crystals) are
defined by axioms (A1)--(A5) below; we give axiomatics in a slightly
different, but equivalent, form compared with~\cite{ste-03}. In what
follows an $n$-colored digraph $K$ is assumed to be a (weakly) {\em
connected}, i.e., it is not representable as the disjoint union of two
nonempty digraphs.

The first axiom concerns the structure of monochromatic subgraphs
$(V,E_i)$.
\begin{itemize}
\item[(A1)] For $i=1,\ldots,n$, each maximal connected subgraph (component)
of $(V,E_i)$ is a simple {\em finite} path, i.e., a sequence of
the form $(v_0,e_1,v_1,\ldots,e_k,v_k)$, where
$v_0,v_1,\ldots,v_k$ are distinct vertices and each $e_i$ is an edge from
$v_{i-1}$ to $v_i$.
  \end{itemize}

In particular, for each $i$, each vertex has at most one incoming $i$-edge
and at most one outgoing $i$-edge, and therefore, one can associate to the
set $E_i$ partial invertible operator $F_i$ acting on vertices: $(u,v)$ is
an $i$-edge if and only if $F_i$ is applicable to $u$ and $F_i(u)=v$.
Since $K$ is connected, one can use the operator notation to express any
vertex via another one. For example, the expression $F_1^{-1}F_3^2F_2(v)$
(where $F_p^{-1}$ stands for the partial operator inverse to $F_p$)
determines the vertex $w$ obtained from a vertex $v$ by traversing 2-edge
$(v,v')$, followed by traversing 3-edges $(v',u)$ and $(u,u')$, followed
by traversing 1-edge $(w,u')$ in backward direction. Emphasize that every
time we use such an operator expression in what follows, this
automatically indicates that all involved edges do exist in $K$.

We refer to a {\em maximal} monochromatic path with color $i$ on the edges
as an $i$-{\em line}. The $i$-line passing a vertex $v$ (possibly
consisting of the only vertex $v$) is denoted by $P_i(v)$, its part from
the first vertex to $v$ by $\Pin_i(v)$, and its part from $v$ to the last
vertex by $\Pout_i(v)$. The lengths of $\Pin_i(v)$ and $\Pout_i(v)$ (i.e.,
the numbers of edges in these paths) are denoted by $\ellin_i(v)$ and
$\ellout_i(v)$, respectively.

Axioms (A2)--(A5) tell us about interrelations of different colors $i,j$.
Taken together, they are equivalent to saying that each component
of the digraph $(V(K),E_i\cup E_j)$ forms a regular $A_2$-crystal when
colors $i,j$ are {\em neighboring}, i.e., $|i-j|=1$, and forms a regular
$A_1\times A_1$-crystal (the Cartesian product of two paths) otherwise.

The second axiom indicates possible changes of the head and tail part
lengths of $j$-lines when one traverses an edge of another color $i$;
these changes depend on the Cartan matrix.

 \begin{itemize}
 \item[(A2)]
For any two colors $i\ne j$ and for any edge $(u,v)$ with the color $i$,
one holds $\ellin_j(v)\le\ellin_j(u)$ and $\ellout_j(v)\ge\ellout_j(u)$.
The value $(\ellin_j(v)-\ellin_j(u))+(\ellout_j(u)-\ellout_j(v))$ is equal
to the coefficient $m_{ij}$ in the Cartan matrix $M$. Furthermore, $h_j$
is convex on each $i$-path, in the sense that if $(u,v),(v,w)$ are
consecutive $i$-edges, then $h_j(u)+h_j(w)\ge 2h_j(v)$.
 \end{itemize}

This can be rewritten as follows.
  \begin{numitem}
When $|i-j|=1$, each $i$-line $P$ contains a vertex $r$ such that: for any
edge $(u,v)$ in $\Pin_i(r)$, one holds $\ellin_j(v)=\ellin_j(u)-1$ and
$\ellout_j(v)=\ellout_j(u)$, and for any edge $(u',v')$ in $\Pout_i(r)$,
one holds $\ellin_j(v')=\ellin_j(u')$ and $\ellout_j(v')=\ellout_j(u')+1$.
When $|i-j|\ge 2$, any $i$-edge $(u,v)$ satisfies
$\ellin_j(v)=\ellin_j(u)$ and $\ellout_j(v)=\ellout_j(u)$.
  \label{eq:A2conv}
  \end{numitem}

Such a vertex $r$ (which is unique) is called the {\em critical} vertex
for $P,i,j$. It is convenient to assign to each $i$-edge $e$ {\em label}
$\ell_j(e)$ taking value 0 if $e$ occurs in the corresponding $i$-line
{\em before} the critical vertex, and 1 otherwise. Emphasize that the
critical vertex (and therefore, edge labels) on an $i$-line $P$ depends on
$j$: the critical vertices on $P$ with respect to the neighboring colors
$j=i-1$ and $j=i+1$ may be different.

Two operators $F=F_i^{\alpha}$ and $F'=F_j^\beta$, where $\alpha,\beta\in
\{1,-1\}$, are said to {\em commute} at a vertex $v$ if each of $F,F'$ acts at
$v$ and $FF'(v)=F'F(v)$. The third axiom points out the situations when
operators commute for neighboring colors $i,j$.

\begin{itemize}
\item[(A3)] Let $|i-j|=1$. (a) If a vertex $u$ has outgoing $i$-edge $(u,v)$
and outgoing $j$-edge $(u,v')$ and if $\ell_j(u,v)=0$, then $\ell_{i}(u,v')=1$
and $F_iF_j(u)=F_jF_i(u)$. Symmetrically: (b) if a vertex $v$ has incoming
$i$-edge $(u,v)$ and incoming $j$-edge $(u',v)$ and if $\ell_j(u,v)=1$, then
$\ell_{i}(u',v)=0$ and $F_i^{-1}F_j^{-1}(v)=F_j^{-1}F_i^{-1}(v)$. (See the
picture.)
  \end{itemize}
 \begin{center}
  \unitlength=1mm
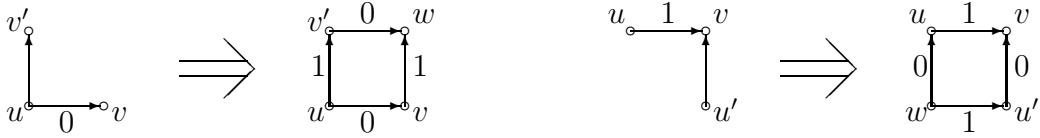
  \begin{picture}(140,20)
\put(5,5){\circle{1.0}}
\put(15,5){\circle{1.0}}
\put(45,5){\circle{1.0}}
\put(55,5){\circle{1.0}}
\put(95,5){\circle{1.0}}
\put(125,5){\circle{1.0}}
\put(135,5){\circle{1.0}}
\put(5,15){\circle{1.0}}
\put(45,15){\circle{1.0}}
\put(55,15){\circle{1.0}}
\put(85,15){\circle{1.0}}
\put(95,15){\circle{1.0}}
\put(125,15){\circle{1.0}}
\put(135,15){\circle{1.0}}
\put(5,5){\vector(1,0){9.5}}
\put(45,5){\vector(1,0){9.5}}
\put(125,5){\vector(1,0){9.5}}
\put(45,15){\vector(1,0){9.5}}
\put(85,15){\vector(1,0){9.5}}
\put(125,15){\vector(1,0){9.5}}
\put(5,5){\vector(0,1){9.5}}
\put(45,5){\vector(0,1){9.5}}
\put(55,5){\vector(0,1){9.5}}
\put(95,5){\vector(0,1){9.5}}
\put(125,5){\vector(0,1){9.5}}
\put(135,5){\vector(0,1){9.5}}
\put(25,9){\line(1,0){9}}
\put(25,11){\line(1,0){9}}
\put(105,9){\line(1,0){9}}
\put(105,11){\line(1,0){9}}
\put(31,6){\line(1,1){4}}
\put(31,14){\line(1,-1){4}}
\put(111,6){\line(1,1){4}}
\put(111,14){\line(1,-1){4}}
 \put(2,3){$u$}
 \put(2,15){$v'$}
 \put(16,3){$v$}
 \put(42,3){$u$}
 \put(42,15){$v'$}
 \put(56,3){$v$}
 \put(56,16){$w$}
 \put(82,16){$u$}
 \put(96,16){$v$}
 \put(96,3){$u'$}
 \put(121.5,3){$w$}
 \put(122,16){$u$}
 \put(136,16){$v$}
 \put(136,3){$u'$}
 \put(9,1.5){0}
 \put(49,1.5){0}
 \put(129,1.5){1}
 \put(49,16){0}
 \put(89,16){1}
 \put(129,16){1}
 \put(42.5,9){1}
 \put(56,9){1}
 \put(122.5,9){0}
\put(136,9){0}
  \end{picture}
 \end{center}

Note that for each ``square'' $u,v,v',w$, where $v=F_i(u)$, $v'=F_j(u)$
and $w=F_j(v)=F_i(v')$, the trivial relations $h_j(u)=h_j(v')+1$ and
$h_j(v)=h_j(w)+1$ imply that the opposite $i$-edges $(u,v)$ and $(v',w)$
have equal labels $\ell_j$; similarly $\ell_{i}(u,v')=\ell_{i}(v,w)$.
Another important consequence of (A3) is that
  \begin{numitem}
for $|i-j|=1$, if $v$ is the critical vertex on an $i$-line with respect
to the color $j$, then $v$ is the critical vertex on the $j$-line passing
$v$ with respect to the color $i$,
  \label{eq:equal_crit}
  \end{numitem}
i.e., we can speak of common critical vertices for the pair $\{i,j\}$.
Indeed, if a vertex $v$ has incoming $i$-edge $(u,v)$ with $\ell_j(u,v)=0$
and outgoing $j$-edge $(v,w)$, then $h_j(u)=h_j(v)\ge 1$, and hence $u$
has outgoing $j$-edge $(u,v')$. By (A3), $w=F_i(v')$ and
$\ell_{i}(u,v')=1$; the latter implies $\ell_{i}(v,w)=1$. Symmetrically,
if $v$ has outgoing $i$-edge $e$ with $\ell_j(e)=1$ and incoming $j$-edge
$e'$, then $\ell_{i}(e')=0$.

The fourth axiom points out the situations when for neighboring $i,j$, the
operators $F_i,F_j$ and their inverse ones  ``remotely commute'' (they are said
to satisfy the ``Verma relation of degree 4'').

\begin{itemize}
\item[(A4)] Let $|i-j|=1$.
(i) If a vertex $u$ has outgoing edges with the colors $i$ and $j$ and if
each edge is labeled 1 (with respect to the other color), then
$F_iF_j^2F_i(u)=F_jF_i^2F_j(u)$. Symmetrically: (ii) if $v$ has incoming
edges with the color $i$ and $j$ and if both are labeled 0, then
$F_i^{-1}(F_j^{-1})^2 F_i^{-1}(v)=F_j^{-1}(F_i^{-1})^2 F_j^{-1}(v)$. (See
the picture.)
  \end{itemize}
 \begin{center}
  \unitlength=1mm
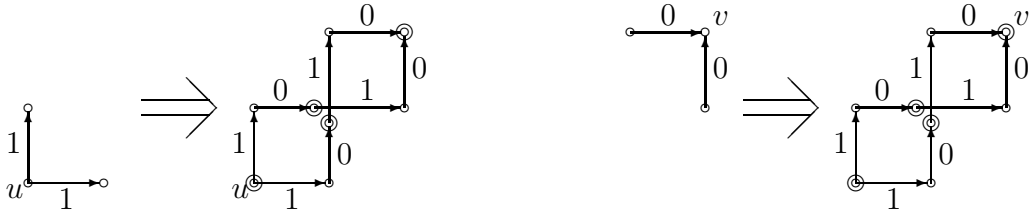
  \begin{picture}(147,30)
\put(5,5){\circle{1.0}}
\put(15,5){\circle{1.0}}
\put(5,15){\circle{1.0}}
\put(2,3){$u$}
\put(5,5){\vector(1,0){9.5}}
\put(5,5){\vector(0,1){9.5}}
\put(9,1.5){1}
\put(2,9){1}
\put(20,14){\line(1,0){9}}
\put(20,16){\line(1,0){9}}
\put(26,11){\line(1,1){4}}
\put(26,19){\line(1,-1){4}}
\put(35,5){\circle{1.0}}
\put(45,5){\circle{1.0}}
\put(35,15){\circle{1.0}}
\put(45,13){\circle{1.0}}
\put(43,15){\circle{1.0}}
\put(55,15){\circle{1.0}}
\put(45,25){\circle{1.0}}
\put(55,25){\circle{1.0}}
\put(35,5){\circle{2.0}}
\put(45,13){\circle{2.0}}
\put(43,15){\circle{2.0}}
\put(55,25){\circle{2.0}}
\put(35,5){\vector(1,0){9.5}}
\put(35,15){\vector(1,0){7.5}}
\put(43,15){\vector(1,0){11.5}}
\put(45,25){\vector(1,0){9.5}}
\put(35,5){\vector(0,1){9.5}}
\put(45,5){\vector(0,1){7.5}}
\put(45,13){\vector(0,1){11.5}}
\put(55,15){\vector(0,1){9.5}}
\put(32,3){$u$}
\put(39,1.5){1}
\put(32,9){1}
\put(46,7.5){0}
\put(37.5,16){0}
\put(49,16){1}
\put(42,19){1}
\put(56,19){0}
\put(49,26){0}
\put(85,25){\circle{1.0}}
\put(95,15){\circle{1.0}}
\put(95,25){\circle{1.0}}
\put(85,25){\vector(1,0){9.5}}
\put(95,15){\vector(0,1){9.5}}
\put(96,26){$v$}
\put(89,26){0}
\put(96,19){0}
\put(100,14){\line(1,0){9}}
\put(100,16){\line(1,0){9}}
\put(106,11){\line(1,1){4}}
\put(106,19){\line(1,-1){4}}
\put(115,5){\circle{1.0}}
\put(125,5){\circle{1.0}}
\put(115,15){\circle{1.0}}
\put(125,13){\circle{1.0}}
\put(123,15){\circle{1.0}}
\put(135,15){\circle{1.0}}
\put(125,25){\circle{1.0}}
\put(135,25){\circle{1.0}}
\put(115,5){\circle{2.0}}
\put(125,13){\circle{2.0}}
\put(123,15){\circle{2.0}}
\put(135,25){\circle{2.0}}
\put(115,5){\vector(1,0){9.5}}
\put(115,15){\vector(1,0){7.5}}
\put(123,15){\vector(1,0){11.5}}
\put(125,25){\vector(1,0){9.5}}
\put(115,5){\vector(0,1){9.5}}
\put(125,5){\vector(0,1){7.5}}
\put(125,13){\vector(0,1){11.5}}
\put(135,15){\vector(0,1){9.5}}
\put(136,26){$v$}
\put(119,1.5){1}
\put(112,9){1}
\put(126,7.5){0}
\put(117.5,16){0}
\put(129,16){1}
\put(122,19){1}
\put(136,19){0}
\put(129,26){0}
  \end{picture}
 \end{center}

One can show that the labels with respect to $i$ or $j$ of all involved
edges are determined uniquely, just as indicated in the above picture
(where the circles indicate the critical vertices).

The final axiom concerns non-neighboring colors.

 \begin{itemize}
 \item[(A5)]
Let $|i-j|\ge 2$. Then for any $F\in\{F_i,F_i^{-1}\}$ and $F'\in
\{F_j,F_j^{-1}\}$, the operators $F,F'$ commute at each vertex where both
act.
  \end{itemize}

This is equivalent to saying that for $|i-j|\ge 2$, each component of the
2-colored subgraph $(V(K),E_i\cup E_j)$ is the Cartesian product of a path
with the color $i$ and a path with the color $j$, i.e., it is an
$A_1\times A_1$-{\em crystal} (viewed as a rectangular grid).

 \bigskip
 \noindent
 {\bf\large \arabic{section}.2. Some properties of RAN-crystals.}
We review some known properties of RAN-crystals that will be used later.

We say that a vertex $v$ of a finite or infinite digraph $G$ is the {\em
source} (resp. {\em sink}) if any inclusion-wise maximal path begins
(resp. ends) at $v$; in particular, $v$ has zero indegree (resp. zero
outdegree). When such a vertex exists, we say that $G$ {\em has source}
(resp. {\em has sink}). The importance of simply-laced crystals with
source is emphasized by a result of Stembridge in~\cite{ste-03}; in the
$A_n$ case it reads as follows:
  \begin{numitem}
For any $n$-tuple $\bfc=(c_1,\ldots,c_n)$ of nonnegative integers, there
exists precisely one RAN-crystal $K$ with source $s$ such that
$\ellout_i(s)=c_i$ for $i=1,\ldots,n$. This $K$ is the crystal graph of
the integrable $U_q(sl_{n+1})$-module of highest weight $\bfc$.
  \label{eq:sourceAn}
  \end{numitem}
(Hereinafter we usually denote $n$-tuples in bold.) We say that $\bfc$ is the
{\em parameter} (tuple) of such a $K$ and denote $K$ by $K(\bfc)$. If we
reverse the edges of $K$ while preserving their colors, we again obtain an
RAN-crystal (since (A1)--(A5) remain valid for it). It is called the {\em dual}
of $K$ and denoted by $K^{\ast}$.

Another property, indicated in~\cite{ste-03} for simply-laced crystals
with a nonsingular Cartan matrix, is easy.
  \begin{numitem}
An RAN-crystal $K$ is {\em graded} for each color $i$, which means that
for any cycle ignoring the orientation of edges, the number of $i$-edges
in one direction is equal to the number of $i$-edges in the other
direction. (One also says that $K$ admits a weight mapping.) In
particular, $K$ is acyclic and has no parallel edges.
  \label{eq:grad}
  \end{numitem}
(Indeed, associate to each vertex $v$ the $n$-vector $wt(v)$
whose $j$-th entry is equal to $\ellout_j(v)-\ellin_j(v)$,
$j=1,\ldots,n$. Then for each $i$-edge $(u,v)$, the difference
$wt(u)-wt(v)$ coincides with the $i$-th row vector $m_i$ of
the Cartan matrix $M$, in view of axiom (A2) and the obvious equality
$\ellout_i(u)-\ellin_i(u)= \ellout_i(v)-\ellin_i(v)+2$. So under the
map $wt:V(K)\to\Rset^n$, the edges of each color $i$ correspond
to parallel translations of one and the same vector $m_i$, and
now~\refeq{grad} follows from the fact that the vectors
$m_1,\ldots,m_n$ are linearly independent.)

\medskip
In general a regular simply-laced crystal need not have source and/or
sink; it may be infinite and may contain directed cycles. One
simple result on regular simply-laced crystals in~\cite{ste-03}
remains valid for more general digraphs, in
particular, for a larger class of crystals of representations.
  \begin{prop}  \label{pr:cr_source}
Let $G$ be an (uncolored) connected and graded digraph with the
following property ($\ast$): for any vertex $v$ and any edges $e,e'$
entering $v$, there exist two paths from some vertex $w$ to $v$ such
that one path contains $e$ and the other contains $e'$. Then either
$G$ has source or all maximal paths in $G$ are infinite in backward
direction.
   \end{prop}
(A similar assertion concerns sinks and infinite paths in forward
direction. For any RAN-crystal, condition ($\ast$) in the proposition is
provided by axioms (A3)--(A5).)

\smallskip
   \begin{proof}
Suppose this is not so. Then, since $G$ is connected and acyclic (as it is
graded), there exists a vertex $v$ and two paths $P,P'$ ending at $v$ such
that $P$ begins at a zero-indegree vertex $s$, while $P'$ either is
infinite in backward direction or begins at a zero-indegree vertex
different from $s$. Let such $v,P,P'$ be chosen so that the length $|P|$
of $P$ is minimum. Then the last edges $e=(u,v)$ and $e'=(u',v)$ of $P$
and $P'$, respectively, are different. By ($\ast$), there is a vertex $w$,
a path $Q$ from $w$ to $v$ containing $e$ and a path $Q'$ from $w$ to $v$
containing $e'$. Extend $Q$ to a maximal path $\bar Q$ ending at $v$.
Three cases are possible: (i) $\bar Q$ is infinite in backward direction;
(ii) $\bar Q$ begins at a (zero-indegree) vertex different from $s$; and
(iii) $\bar Q$ begins at $s$. In cases~(i),(ii), we come to a
contradiction with the minimality of $P$ by taking the vertex $u$ and the
part of $P$ from $s$ to $u$. And in case~(iii), there is a path $\bar Q'$
from $s$ to $v$ that contains $e'$. Since $G$ is graded, $|\bar Q'|=|P|$.
Then we again get a contradiction with the minimality of $P$ by taking
$u'$, the part of $\bar Q'$ from $s$ to $u'$, and the part of $P'$ ending
at $u'$.
   \end{proof}

(The fact that $G$ is graded is important. Indeed, take $G$ with the
vertices $s$ and $u_i,v_i$ for all $i\in\Zset_+$, and the edges $(s,u_0)$
and $(u_i,u_{i+1}),(v_{i+1},v_i),(u_i,v_i)$ for all $i$. This $G$
satisfies ($\ast$), the vertex $s$ has zero indegree, and the path on the
vertices $v_i$ is infinite in backward direction. One can also construct a
locally finite graph satisfying ($\ast$) and having many zero-indegree
vertices.)

Our crossing model will generate $n$-colored graphs satisfying axioms
(A1)--(A5); moreover, it generates one RAN-crystal with source for each
parameter tuple $\bfc\in\Zset_+^n$. In light of~\refeq{sourceAn} and
Proposition~\ref{pr:cr_source}, a reasonable question is whether {\em
every} RAN-crystal has source and sink (or, equivalently, is finite). The
question will be answered affirmatively in the next subsection, thus
implying that the crossing model gives the whole set of RAN-crystals.

As a consequence of the crossing model, we also will observe the following
anti-symmetric property of an RAN-crystal $K$: if we reverse the numeration of
colors (regarding each color $i$ as $n-i+1$) in the dual crystal $K^\ast$, then
the resulting crystal is isomorphic to $K$. In other words, $\ellout_i(s_K)=
\ellin_{n-i+1}(\bar s_K)$ for $i=1,\ldots,n$, where $s_K$ and $\bar s_K$ are
the source and sink of $K$, respectively.

\medskip
Finally, recall that a Gelfand-Tsetlin pattern~\cite{GT-50}, or a {\em
GT-pattern} for short, is a triangular array $X=(x_{ij})_{1\le j\le i\le
n}$ of integers satisfying $x_{ij}\ge x_{i-1,j},x_{i+1,j+1}$ for all
$i,j$. Given a weakly decreasing $n$-tuple $\bfa=(a_1\ge \cdots\ge a_n)$
of nonnegative integers, one says that $X$ {\em is bounded} by $\bfa$ if
$a_j\ge x_{n,j}\ge a_{j+1}$ for $j=1,\ldots,n$, letting $a_{n+1}:=0$. It
is known that GT-patterns, as well as the corresponding semi-standard
Young tableaux, are closely related to crystals of representations for
$U_q(sl_{n+1})$ (cf.~\cite{BZ-01,kas-95,KN-94,lit-95b}). More precisely,
  \begin{numitem}
for any $\bfc\in\Zset_+^n$, there is a bijection between the vertex set of
the RAN-crystal $K(\bfc)$ and the set of GT-patterns bounded by the
$n$-tuple $\bfc^\Sigma=(c^\Sigma_1,\ldots,c^\Sigma_n)$, defined by
$c^\Sigma_j:=c_1+\ldots+c_{n-j+1}$ for $j=1,\ldots,n$.
  \label{eq:cr-GT}
  \end{numitem}
As mentioned in the Introduction, there is a correspondence between
GT-patterns and feasible functions in the crossing model; it will be
exposed in Proposition~\ref{pr:f-patt}.

 \bigskip
 \noindent
 {\bf\large \arabic{section}.3. Properties of $A_2$-crystals.}
In this subsection we give a brief review of certain results
from~\cite{DKK-07} for the simplest case $n=2$, namely, for regular
$A_2$-crystals, or {\em RA2-crystals} for short. They describe the
combinatorial structure of such crystals and demonstrate some additional
properties.

An RA2-crystal $K$ is defined by axioms (A1)--(A4) with $\{i,j\}=\{1,2\}$
(since (A5) becomes redundant). It turns out that these crystals can be
produced from elementary 2-colored crystals by use of a certain operation
of replicating and gluing together. This operation can be introduced for a
pair of arbitrary finite or infinite graphs as follows. (In
Section~\ref{sec:princ} the construction is generalized to $n$ graphs, in
connection with the so-called skeleton of an RAN-crystal.)

Consider graphs $G=(V,E)$ and $H=(V',E')$ with distinguished vertex
subsets $S\subseteq V$ and $T\subseteq V'$. Take $|T|$ disjoint copies of
$G$, denoted as $G_t$ ($t\in T$), and $|S|$ disjoint copies of $H$,
denoted as $H_s$ ($s\in S$). We glue these copies together in the
following way: for each $s\in S$ and each $t\in T$, the vertex $s$ in
$G_t$ is identified with the vertex $t$ in $H_s$. The resulting graph,
consisting of $|V| |T|+|V'| |S|-|S| |T|$ vertices and $|E| |T|+|E'| |S|$
edges, is denoted as $(G,S)$\sprod$(H,T)$.

In our special case the role of $G$ and $H$ is played by 2-colored
digraphs $R$ and $L$ viewed as triangular parts of square grids. More
precisely, $R$ depends on a parameter $c_1\in\Zset_+$ and its vertices $v$
correspond to the integer points $(i,j)$ in the plane such that $0\le j\le
i\le c_1$. The vertices $v$ of $L$, depending on a parameter
$c_2\in\Zset_+$, correspond to the integer points $(i,j)$ such that $0\le
i\le j\le c_2$. We say that $v$ has the {\em coordinates} $(i,j)$ in the
sail. The edges with the color 1 in these digraphs correspond to all
possible pairs $((i,j),(i+1,j))$, and the edges with the color 2 to the
pairs $((i,j),(i,j+1))$. We call $R$ the {\em right sail} of size $c_1$,
and $L$ the {\em left sail} of size $c_2$.

It is easy to check
that $R$ satisfies axioms (A1)--(A4) and is just the crystal $K(c_1,0)$,
and that the set of critical vertices in $R$ coincides with the {\em diagonal}
$D_R=\{(i,i): i=0,\ldots,c_1\}$. Similarly, $L=K(0,c_2)$,
and the set of critical vertices in it coincides with the diagonal
$D_L=\{(i,i): i=0,\ldots,c_2\}$. These diagonals are just taken as the
distinguished subsets in these digraphs. The vertices in
$D_R$ ($D_L$) are ordered in a natural way, according to which $(i,i)$ is
referred as the $i$-th critical vertex in $R$ ($L$).

We refer to the digraph obtained by use of operation \sprod in this case
as the {\em diagonal-product} of $R$ and $L$, and for brevity write
$R$\sprod$L$, omitting the distinguished subsets. The edge colors in the
resulting graph are inherited from $R$ and $L$. Using the above ordering
in the diagonals, we may speak of $p$-th right sail in $R$\sprod$L$,
denoted by $R_p$. Here $0\le p\le c_2$, and $R_p$ is the copy of $R$
corresponding to the vertex $(p,p)$ of $L$. In a similar way, one defines
$q$-th left sail $L_q$ in $R$\sprod$L$ for $q=0,\ldots,c_1$. The common
vertex of $R_p$ and $L_q$ is denoted by $v_{p,q}$.

One checks that $R$\sprod$L$ has source and sink and satisfies axioms
(A1)--(A4). Moreover, it is exactly the RA2-crystal $K(c_1,c_2)$. The
critical vertices in it are just $v_{p,q}$ for all $p,q$, the source is
$v_{0,0}$ and the sink is $v_{c_1,c_2}$. The case $c_1=1$ and $c_2=2$ is
illustrated in Fig.~\ref{fig:sprod}; here the critical vertices are
indicated by circles, 1-edges by horizontal arrows, and 2-edges by
vertical arrows.

\begin{figure}[htb]                
 \begin{center}
  \unitlength=1mm
  \begin{picture}(125,35)
\put(5,5){\circle{1.0}}           
\put(5,15){\circle{1.0}}
\put(5,25){\circle{1.0}}
\put(15,15){\circle{1.0}}
\put(15,25){\circle{1.0}}
\put(25,25){\circle{1.0}}
\put(5,5){\circle{2.5}}
\put(15,15){\circle{2.5}}
\put(25,25){\circle{2.5}}
\put(5,15){\vector(1,0){9.5}}
\put(5,25){\vector(1,0){9.5}}
\put(15,25){\vector(1,0){9.5}}
\put(5,5){\vector(0,1){9.5}}
\put(5,15){\vector(0,1){9.5}}
\put(15,15){\vector(0,1){9.5}}
\put(20,5){(a)}
\put(55,5){\circle{1.0}}
\put(65,5){\circle{1.0}}
\put(65,15){\circle{1.0}}
\put(55,5){\circle{2.5}}
\put(65,15){\circle{2.5}}
\put(55,5){\vector(1,0){9.5}}
\put(65,5){\vector(0,1){9.5}}
\put(70,5){(b)}
\put(95,5){\circle{1.0}}
\put(95,15){\circle{1.0}}
\put(95,25){\circle{1.0}}
\put(105,15){\circle{1.0}}
\put(105,25){\circle{1.0}}
\put(115,25){\circle{1.0}}
\put(95,5){\circle{2.5}}
\put(105,15){\circle{2.5}}
\put(115,25){\circle{2.5}}
\put(95,15){\vector(1,0){9.5}}
\put(95,25){\vector(1,0){9.5}}
\put(105,25){\vector(1,0){9.5}}
\put(95,5){\vector(0,1){9.5}}
\put(95,15){\vector(0,1){9.5}}
\put(105,15){\vector(0,1){9.5}}
\put(102,12){\circle{1.0}}
\put(102,22){\circle{1.0}}
\put(102,32){\circle{1.0}}
\put(112,22){\circle{1.0}}
\put(112,32){\circle{1.0}}
\put(122,32){\circle{1.0}}
\put(102,12){\circle{2.5}}
\put(112,22){\circle{2.5}}
\put(122,32){\circle{2.5}}
\put(102,22){\vector(1,0){9.5}}
\put(102,32){\vector(1,0){9.5}}
\put(112,32){\vector(1,0){9.5}}
\put(102,12){\vector(0,1){9.5}}
\put(102,22){\vector(0,1){9.5}}
\put(112,22){\vector(0,1){9.5}}
\put(102,5){\circle{1.0}}
\put(95,5){\vector(1,0){6.5}}
\put(102,5){\vector(0,1){6.5}}
\put(112,15){\circle{1.0}}
\put(105,15){\vector(1,0){6.5}}
\put(112,15){\vector(0,1){6.5}}
\put(122,25){\circle{1.0}}
\put(115,25){\vector(1,0){6.5}}
\put(122,25){\vector(0,1){6.5}}
\put(115,5){(c)}
  \end{picture}
 \end{center}
 \caption{(a) $K(0,2)$,\;\; (b) $K(1,0)$,\;\;
(c) $K(1,0)\triangleright\hspace{-4pt}\triangleleft\, K(0,2)$.}
  \label{fig:sprod}
  \end{figure}

\begin{theorem} \label{tm:main2} {\rm \cite{DKK-07}}
Any RA2-crystal $K$ is representable as $K(a,0)$\sprod$K(0,b)$ for some
$a,b\in\Zset_+$ (in particular, $K$ is finite). The set of RA2-crystals is
exactly $\{K(\bfc): \bfc\in \Zset_+^2\}$.
  \end{theorem}

A useful consequence of the above construction is that the vertices $v$ of
$K$ one-to-one correspond to the quadruples $(\alpha_1,\alpha_2,\beta_1,
\beta_2)$ of integers such that
   \begin{numitem}
(i) $0\le \alpha_2\le\alpha_1\le c_1$, (ii) $0\le\beta_1\le\beta_2\le c_2$,
and (iii) at least one of the equalities $\alpha_2=\alpha_1$ and
$\beta_1=\beta_2$ takes place,
   \label{eq:alpha-beta}
  \end{numitem}
and each i-edge ($i=1,2$) corresponds to the increase by 1 of one of
$\alpha_i,\beta_i$, subject to maintaining~\refeq{alpha-beta}.

Under this correspondence, if $\beta_1=\beta_2$ then $v$ occurs in the
right sail with the number $\beta_1$ and has the coordinates
$(\alpha_1,\alpha_2)$ in it, while if $\alpha_2=\alpha_1$ then $v$ occurs
in the left sail with the number $\alpha_1$ and has the coordinates
$(\beta_1,\beta_2)$. In particular, a critical vertex $v_{p,q}$
corresponds to $(q,q,p,p)$.

\medskip
\noindent {\bf Remark 1.} The representation of the vertices of $K$ as the
above quadruples satisfying~\refeq{alpha-beta} gives rise to constructing
the crossing model for the simplest case $n=2$, as we explain in the next
section. A more general numerical representation (which is beyond our
consideration in this paper) does not impose condition (iii)
in~\refeq{alpha-beta}. In this case the admissible transformations of
quadruples $(\alpha_1,\alpha_2,\beta_1,\beta_2)$ (giving the edges of a
digraph on the quadruples) are assigned as follows. For $\Delta:=
\min\{\alpha_1-\alpha_2,\beta_2-\beta_1\}$, we choose one of
$\alpha_1,\alpha_2,\beta_1,\beta_2$ and increase it by 1 unless this
increase violates (i) or (ii) in~\refeq{alpha-beta} or changes $\Delta$.
One can see that the resulting digraph $Q$ is the disjoint union of
$1+\min\{c_1,c_2\}$ RA2-crystals, namely, $K(c_1-\Delta,c_2-\Delta)$ for
$\Delta=0,\ldots,\min\{c_1,c_2\}$. (This $Q$ is the tensor product of
crystals (sails) $K(c_1,0)$ and $K(0,c_2)$.)

\medskip
One more useful result in~\cite{DKK-07} is as follows.
  \begin{prop} \label{pr:A4}
Part (ii) of axiom (A4) for RAN-crystals is redundant. Furthermore, axiom
(A4) itself follows from (A1)--(A3) if we add the condition that each
component of $(V,E_i\cup E_j)$ with $|i-j|=1$ has exactly one
zero-indegree (or exactly one zero-outdegree) vertex.
  \end{prop}

 \medskip
In conclusion of this section, return to an arbitrary RAN-crystal $K$. For
a color $i$, let $H_i$ denote the operator on $V(K)$ that brings a vertex
$v$ to the end vertex of the path $P_i(v)$, i.e., $H_i(v)=F_i^{h_i(v)}(v)$
(letting $F_i^0={\rm id}$). We observe that
  \begin{numitem}
for neighboring colors $i,j$ and a vertex $v$, if $h_i(v)=0$ then the
vertex $w=H_iH_j(v)$ satisfies $h_i(w)=h_j(w)=0$.
   \label{eq:HiHj}
  \end{numitem}

Indeed, the RA2-subcrystal with the colors $i,j$ in $K$ that contains $v$
is $K(c_i,c_j)$ for some $c_i,c_j\in\Zset_+$. Represent $v$ as quadruple
$q=(\alpha_i,\alpha_j,\beta_i,\beta_j)$ in~\refeq{alpha-beta} (with $i,j$
in place of 1,2). Then $h_i(q)=0$ implies $\alpha_i=c_i$ and
$\beta_i=\beta_j$. One can see that applying $H_j$ to $q$ results in the
quadruple $q'=(c_i,c_i,\beta_i,c_j)$ and applying $H_i$ to $q'$ results in
$(c_i,c_i,c_j,c_j)$. This gives~\refeq{HiHj}.

Using~\refeq{HiHj}, we can show the following important property of
RAN-crystals.
  \begin{prop}  \label{pr:RANsink}
Any RAN-crystal $K$ has a zero-outdegree vertex.
  \end{prop}
  \begin{proof}
For a vertex $u$, let $p(u)$ be the maximum integer $p$ such that
$h_i(u)=0$ for $i=1,\ldots,p-1$. Assuming $p(u)<n+1$, we claim that the
vertex $w=H_1H_2\ldots H_{p(u)}(u)$ satisfies $p(w)>p(u)$, whence the
result will immediatelly follow. (In other words, by applying the operator
$\bar H_n\bar H_{n-1}\ldots\bar H_1$ to an arbitrary vertex, we get a
zero-outdegree vertex, where $\bar H_i$ stands for $H_1H_2\ldots H_i$.)

Indeed, let $p=p(u)$. For the vertex $v_p:=H_p(u)$, we have $h_p(v_p)=0$
and $h_i(v_p)=h_i(u)$ for all $i\ne p-1,p+1$ (since the colors $p,i$
commute), while $h_{p-1}(v_p)$ may differ from $h_{p-1}(u)$. So
$h_i(v_p)=0$ for $i=1,\ldots,p-2,p$. Similarly, the vertex
$v_{p-1}:=H_{p-1}(v_p)$ satisfies $h_{p-1}(v_{p-1})=0$ and
$h_i(v_{p-1})=h_i(v_p)$ for all $i\ne p-2,p$. Moreover,
applying~\refeq{HiHj} to $v=u$, $i=p-1$ and $j=p$, we obtain
$h_p(v_{p-1})=0$. So $h_i(v_{p-1})=0$ for $i=1,\ldots,p-3,p-1,p$. On the
next step, in a similar fashion one shows that $v_{p-2}:=
H_{p-2}(v_{p-1})$ satisfies $h_i(v_{p-2})=0$ for all $i\in
\{1,\ldots,p\}\setminus\{p-3\}$, and so on. Then the final vertex
$v_1:=H_1\ldots H_p(u)$ in the process has the property $h_i(v_1)=0$ for
$i=1,\ldots,p$, as required in the claim.
  \end{proof}

Also $K$ has a zero-indegree vertex (since
Proposition~\ref{pr:RANsink} can be applied to the dual crystal
$K^\ast$). This together with~\refeq{sourceAn} and
Proposition~\ref{pr:cr_source} gives the following.
  \begin{corollary}  \label{cor:RANsource}
Every RAN-crystal $K$ is finite and has source and sink. Therefore,
$K=K(\bfc)$ for some $\bfc\in\Zset_+^n$.
  \end{corollary}

\section{Description of the crossing model} \label{sec:model}

As mentioned in the Introduction, the {\em crossing model} $\Mscr_n$ for
RAN-crystals consists of three ingredients:

(i) a certain digraph $G=(V(G),E(G)$ depending only on the number $n$ of
colors, called the {\em supporting graph} (the {\em structural} part of
$\Mscr_n$);

(ii) a certain set $\Fscr=\Fscr(\bfc)$ of nonnegative integer-valued
functions on $V(G)$, called {\em feasible functions}, depending on an
$n$-tuple of parameters $\bfc\in\Zset_+^n$ (the {\em numerical} part);

(iii) $n$ partial operators acting on $\Fscr$, called {\em moves} (the
{\em operator} part).

(The feasible functions will correspond to the vertices of the crystal
with the parameter $\bfc$, and the moves to the edges of this crystal.)
Parts (i) and (ii) are described in this section, and part (iii) in the
next one. To avoid a possible mess when both a crystal and the supporting
graph are considered simultaneously, we will refer to a vertex of the
latter graph as a {\em node}.

To explain the idea, we first consider the simplest case $n=2$ and a
2-colored crystal $K=K(c_1,c_2)$. The model $\Mscr_2$ is constructed by
relying on encoding~\refeq{alpha-beta} of the vertices of $K$. The
supporting graph $G$ is formed by two disjoint edges $(u_1,u_2)$ and
$(w_2,w_1)$ (which are related to the elementary crystals, or sails,
$K(c_1,0)$ and $K(0,c_2)$). A feasible function $f$ on $V(G)$ takes values
$f(u_1)=\alpha_1$, $f(u_2)=\alpha_2$, $f(w_1)=\beta_1$, $f(w_2)=\beta_2$
for $\alpha,\beta$ as in~\refeq{alpha-beta}. So the direction of each edge
$e$ of $G$ indicates the corresponding inequality to be imposed on the
values of any feasible function $f$ on the end nodes of $e$, and each $f$
one-to-one corresponds to a vertex of $K$. The graph $G$ is illustrated on
the picture:
 \begin{center}
  \unitlength=1mm
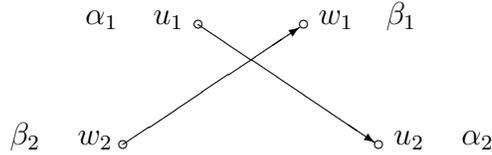
  \begin{picture}(65,18)
\put(15,0){\circle{1.0}}
\put(49,0){\circle{1.0}}
\put(25,16){\circle{1.0}}
\put(39,16){\circle{1.0}}
\put(15,0){\vector(3,2){23.5}}
\put(25,16){\vector(3,-2){23.5}}
\put(9,0){$w_2$}
\put(51,0){$u_2$}
\put(19,16){$u_1$}
\put(41,16){$w_1$}
\put(0,0){$\beta_2$}
\put(60,0){$\alpha_2$}
\put(10,16){$\alpha_1$}
\put(50,16){$\beta_1$}
  \end{picture}
 \end{center}

Note that each admissible quadruple $(\alpha_1,\alpha_2,\beta_1,\beta_2)$
generates the GT-pattern $X$ of size 2 (see Subsection~\ref{sec:An}.3),
defined by $x_{11}:=\alpha_1+\beta_1$, $x_{21}:= \beta_2+c_1$ and
$x_{22}:=\alpha_2$ (see the diagram below). This pattern is bounded by
$\bfc^\Sigma=(c_1+c_2,c_1)$.
 \begin{center}
  \unitlength=1mm
  \begin{picture}(30,12)
\put(0,0){$c_1+\beta_2$}
\put(25,0){$\alpha_2$}
\put(8,10){$\alpha_1+\beta_1$}
  \end{picture}
 \end{center}

Next we start describing the model for an arbitrary $n$. The  ``simplest''
case of an $n$-colored graph $K=K(\bfc)$ arises when all entries in
$\bfc=(c_1,\ldots,c_n)$ are zero except for one entry $c_k$. In this case
we say that $K$ is the $k$-th {\em base crystal} of size $c_k$ and denote
it by $K_n^k(c_k)$.

 \bigskip
 \noindent
 {\bf\large \arabic{section}.1. The supporting graph of $\Mscr_n$.}
To facilitate understanding the construction of the supporting graph $G$,
we first introduce an auxiliary digraph $\Gscr=\Gscr_n$, called the {\em
proto-graph} of $G$. Its node set consists of elements $V_i(j)$ for all
$1\le j\le i\le n$. Its edge set consists of all possible pairs of the
form $(V_i(j),V_{i-1}(j))$ ({\em ascending} edges) or
$(V_i(j),V_{i+1}(j+1))$ ({\em descending} edges). We say that the nodes
$V_i(1),\ldots,V_i(i)$ form the $i$-th {\em level} of $\Gscr$ and order
them as indicated (by increasing $j$). We visualize $\Gscr$ by drawing it
on the plane so that the nodes of the same level lie on a horizontal line,
the edges have equal lengths, the ascending edges point North-East, and
the descending edges point South-East. See the picture for $n=4$.
  \begin{center}
  \unitlength=1mm
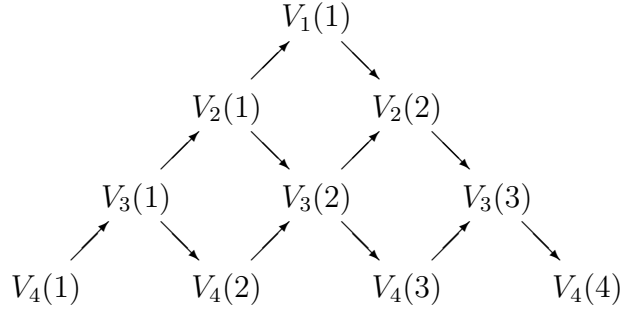
  \begin{picture}(80,40)
   \put(0,0){$V_4(1)$}
   \put(24,0){$V_4(2)$}
   \put(48,0){$V_4(3)$}
   \put(72,0){$V_4(4)$}
   \put(12,12){$V_3(1)$}
   \put(36,12){$V_3(2)$}
   \put(60,12){$V_3(3)$}
   \put(24,24){$V_2(1)$}
   \put(48,24){$V_2(2)$}
   \put(36,36){$V_1(1)$}
  \put(8,5){\vector(1,1){5}}
  \put(32,5){\vector(1,1){5}}
  \put(56,5){\vector(1,1){5}}
  \put(20,17){\vector(1,1){5}}
  \put(44,17){\vector(1,1){5}}
  \put(32,29){\vector(1,1){5}}
  \put(20,10){\vector(1,-1){5}}
  \put(44,10){\vector(1,-1){5}}
  \put(68,10){\vector(1,-1){5}}
  \put(32,22){\vector(1,-1){5}}
  \put(56,22){\vector(1,-1){5}}
  \put(44,34){\vector(1,-1){5}}
  \end{picture}
 \end{center}

The supporting graph $G$ is formed by replicating elements of $\Gscr$ as
follows. Each node $V_i(j)$ generates $n-i+1$ nodes of $G$, denoted as
$v_i^k(j)$ for $k=i-j+1,\ldots, n-j+1$, which are ordered by increasing
$k$ (and accordingly follow from left to right in the visualization). We
identify $V_i(j)$ with the set of these nodes and call it a {\em
multinode} of $G$. Each edge of $\Gscr$ generates a set of edges of $G$ (a
{\em multi-edge}) connecting the elements with equal upper indexes. More
precisely, $(V_i(j),V_{i-1}(j))$ gives $n-i+1$ ascending edges
$(v_i^k(j),v_{i-1}^k(j))$ for $k=i-j+1,\ldots, n-j+1$, and
$(V_i(j),V_{i+1}(j+1))$ gives $n-i$ descending edges
$(v_i^k(j),v_{i+1}^k(j+1))$ for $k=i-j+1,\ldots,n-j$.

The resulting $G$ is the disjoint union of $n$ digraphs $G^1,\ldots,G^n$.
Here $G^k=G^k_n$ contains all vertices of the form $v_i^k(j)$ (the indexes
$i,j$ range over $1\le j\le n-k+1$ and $0\le i-j\le k-1$, and $G^k$ is
viewed as a square (or, better to say, rhombic) grid of size $k-1$ by
$n-k$; we shall see later that $G^k$ is, in fact, the supporting graph for
the base crystal $K^k_n$.) For example: for $n=4$, the graph $G$ is viewed
as
 \begin{center}
  \unitlength=1mm
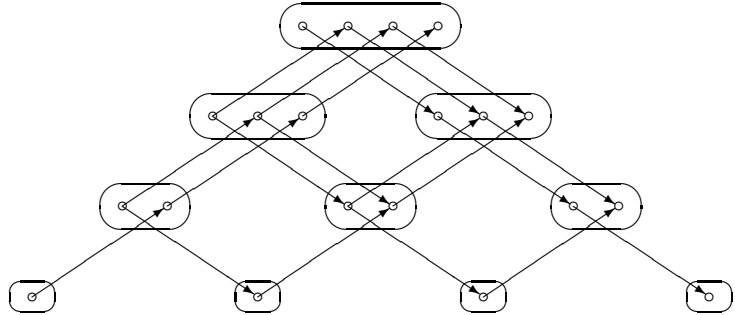
  \begin{picture}(90,39)
\put(0,0){\begin{picture}(54,36)%
\put(0,0){\circle{1.0}} \put(18,12){\circle{1.0}}
\put(36,24){\circle{1.0}} \put(54,36){\circle{1.0}}
\put(0,0){\vector(3,2){17.5}} \put(18,12){\vector(3,2){17.5}}
\put(36,24){\vector(3,2){17.5}}
  \end{picture}}
\put(12,0){\begin{picture}(54,36)%
\put(0,12){\circle{1.0}} \put(18,0){\circle{1.0}}
\put(18,24){\circle{1.0}} \put(36,12){\circle{1.0}}
\put(36,36){\circle{1.0}} \put(54,24){\circle{1.0}}
\put(0,12){\vector(3,2){17.5}} \put(0,12){\vector(3,-2){17.5}}
\put(18,0){\vector(3,2){17.5}} \put(18,24){\vector(3,2){17.5}}
\put(18,24){\vector(3,-2){17.5}} \put(36,12){\vector(3,2){17.5}}
\put(36,36){\vector(3,-2){17.5}}
  \end{picture}}
\put(24,0){\begin{picture}(54,36)%
\put(0,24){\circle{1.0}} \put(18,12){\circle{1.0}}
\put(18,36){\circle{1.0}} \put(36,0){\circle{1.0}}
\put(36,24){\circle{1.0}} \put(54,12){\circle{1.0}}
\put(0,24){\vector(3,2){17.5}} \put(0,24){\vector(3,-2){17.5}}
\put(18,12){\vector(3,2){17.5}} \put(18,12){\vector(3,-2){17.5}}
\put(18,36){\vector(3,-2){17.5}} \put(36,0){\vector(3,2){17.5}}
\put(36,24){\vector(3,-2){17.5}}
  \end{picture}}
\put(36,0){\begin{picture}(54,36)%
\put(0,36){\circle{1.0}} \put(18,24){\circle{1.0}}
\put(36,12){\circle{1.0}} \put(54,0){\circle{1.0}}
\put(0,36){\vector(3,-2){17.5}} \put(18,24){\vector(3,-2){17.5}}
\put(36,12){\vector(3,-2){17.5}}
  \end{picture}}
\put(45,36){\oval(24,6)} \put(30,24){\oval(18,6)} \put(60,24){\oval(18,6)}
\put(15,12){\oval(12,6)} \put(45,12){\oval(12,6)} \put(75,12){\oval(12,6)}
\put(0,0){\oval(6,4)} \put(30,0){\oval(6,4)} \put(60,0){\oval(6,4)}
\put(90,0){\oval(6,4)}
\end{picture}
 \end{center}
 \noindent
(where the multinodes are surrounded by ovals) and its components
$G^1,G^2,G^3,G^4$, called the {\em base subgraphs} of $G$, are viewed as
 \begin{center}
  \unitlength=1mm
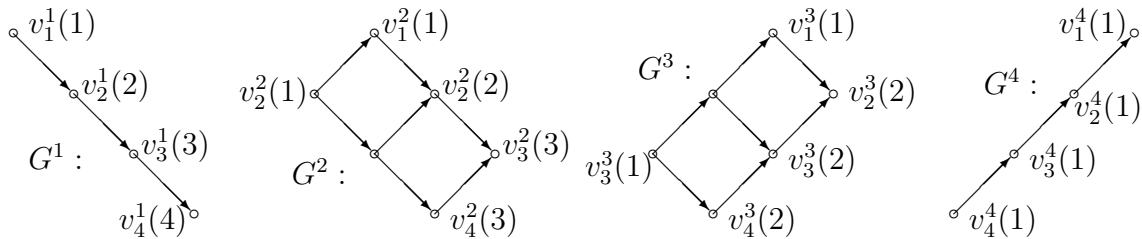
  \begin{picture}(150,27)
\put(0,24){\circle{1.0}} \put(8,16){\circle{1.0}} \put(16,8){\circle{1.0}}
\put(24,0){\circle{1.0}} \put(0,24){\vector(1,-1){7.5}}
\put(8,16){\vector(1,-1){7.5}} \put(16,8){\vector(1,-1){7.5}}
\put(2,6){$G^1:$}
 \put(2,24){$v_1^1(1)$}
 \put(9,16){$v_2^1(2)$}
 \put(17,8){$v_3^1(3)$}
 \put(14,-2){$v_4^1(4)$}
\put(40,16){\circle{1.0}} \put(48,8){\circle{1.0}}
\put(48,24){\circle{1.0}} \put(56,0){\circle{1.0}}
\put(56,16){\circle{1.0}} \put(64,8){\circle{1.0}}
\put(40,16){\vector(1,1){7.5}} \put(40,16){\vector(1,-1){7.5}}
\put(48,8){\vector(1,1){7.5}} \put(48,8){\vector(1,-1){7.5}}
\put(48,24){\vector(1,-1){7.5}} \put(56,0){\vector(1,1){7.5}}
\put(56,16){\vector(1,-1){7.5}} \put(37,4){$G^2:$}
 \put(30,15){$v_2^2(1)$}
 \put(65,8){$v_3^2(3)$}
 \put(57,16){$v_2^2(2)$}
 \put(49,24){$v_1^2(1)$}
 \put(58,-2){$v_4^2(3)$}
\put(85,8){\circle{1.0}} \put(93,0){\circle{1.0}}
\put(93,16){\circle{1.0}} \put(101,8){\circle{1.0}}
\put(101,24){\circle{1.0}} \put(109,16){\circle{1.0}}
\put(85,8){\vector(1,1){7.5}} \put(85,8){\vector(1,-1){7.5}}
\put(93,0){\vector(1,1){7.5}} \put(93,16){\vector(1,1){7.5}}
\put(93,16){\vector(1,-1){7.5}} \put(101,8){\vector(1,1){7.5}}
\put(101,24){\vector(1,-1){7.5}} \put(83,18){$G^3:$}
 \put(76,5){$v_3^3(1)$}
 \put(95,-2){$v_4^3(2)$}
 \put(103,6){$v_3^3(2)$}
 \put(111,15){$v_2^3(2)$}
 \put(103,24){$v_1^3(1)$}

\put(125,0){\circle{1.0}} \put(133,8){\circle{1.0}}
\put(141,16){\circle{1.0}} \put(149,24){\circle{1.0}}
\put(125,0){\vector(1,1){7.5}} \put(133,8){\vector(1,1){7.5}}
\put(141,16){\vector(1,1){7.5}} \put(129,16){$G^4:$}
 \put(127,-2){$v_4^4(1)$}
 \put(135,6){$v_3^4(1)$}
 \put(141,13){$v_2^4(1)$}
 \put(139,24){$v_1^4(1)$}
\end{picture}
 \end{center}

Thus, each node $v=v_i^k(j)$ of $G$ has at most four incident edges,
namely, $(v_{i-1}^k(j-1),v)$, $(v_{i+1}^k(j),v)$, $(v,v_{i-1}^k(j))$,
$(v,v_{i+1}^k(j+1))$; we refer to them, when exist, as the {\em NW-, SW-,
NE-, {\em and} SE-edges}, and denote by $\eNW(v),\eSW(v),\eNE(v),\eSE(v)$,
respectively.

Four nodes of each $G^k$ are distinguished: the leftmost node $v^k_k(1)$,
the rightmost node $v^k_{n-k+1}(n-k+1)$, the topmost node $v^k_1(1)$, and
the bottommost node $v^k_n(n-k+1)$, denoted by
$\lvert^k,~\rvert^k,~\topvert^k$, and $\botvert^k$, respectively. Note
that $\lvert^k$ is the source and $\rvert^k$ is the sink of $G^k$.

 \bigskip
 \noindent
 {\bf\large \arabic{section}.2. Weights of nodes.}
We consider nonnegative integer-valued functions $f$ on $V(G)$ and refer
to the value $f(v)$ as the {\em weight} of a node $v$. A function $f$ is
called {\em feasible} if it satisfies the following three conditions. Here
for an edge $e=(u,v)$, $\;\partial f(e)$ denotes the difference
$f(u)-f(v)$, and $e$ is called {\em tight} for $f$, or $f$-{\em tight}, if
$\partial f(e)=0$.
  \begin{numitem}
  \vspace{-25pt}
 \begin{itemize}
 \item[(i)] $f$ is {\em monotone} on the edges, in the sense that
$\partial f(e)\ge 0$ for all $e\in E(G)$;
 \item[(ii)] $0\le f(v)\le c_k$ for each $v\in V(G^k)$, $k=1,\ldots,n$ (or,
equivalently, $f(\lvert^k)\le c_k$ and $f(\rvert^k)\ge 0$, in view
of~(i));
 \item[(iii)] each multinode $V_i(j)$ contains a node $v$ such that:
the edge $\eSE(u)$ is tight for each node $u\in V_i(j)$ preceding $v$, and
$\eSW(u')$ is tight for each node $u'\in V_i(j)$ succeeding $v$.
  \end{itemize}
 \label{eq:feas}
  \end{numitem}
  \noindent
We say that such a $v$ in (iii) satisfies the {\em switch} condition. The
{\em first} of such nodes $v=v_i^k(j)$ (i.e., with $k$ minimum) is called
the {\em switch-node} in the multinode $V_i(j)$. It plays an important
role in transformations of feasible functions in the model. (We shall see
later that the {\em forward moves}, related to acting operators $F_i$,
handle just switch-nodes, while the {\em backward moves}, related to
acting $F_i^{-1}$, handle {\em last} nodes satisfying the switch
condition.) See the picture, where tight edges are drawn bold and only one
node, marked by a circle, satisfies the switch condition.

 \begin{center}
  \unitlength=1mm
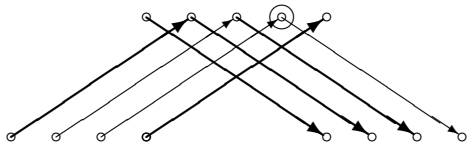
  \begin{picture}(60,18)
 \put(0,0){\circle{1.0}}
 \put(6,0){\circle{1.0}}
 \put(12,0){\circle{1.0}}
 \put(18,0){\circle{1.0}}
 \put(42,0){\circle{1.0}}
 \put(48,0){\circle{1.0}}
 \put(54,0){\circle{1.0}}
 \put(60,0){\circle{1.0}}
 \put(18,0){\circle{1.0}}
 \put(18,16){\circle{1.0}}
 \put(24,16){\circle{1.0}}
 \put(30,16){\circle{1.0}}
 \put(36,16){\circle{1.0}}
 \put(42,16){\circle{1.0}}
 \put(36,16){\circle{3.0}}
 \put(6,0){\vector(3,2){23.5}}
 \put(12,0){\vector(3,2){23.5}}
 \put(36,16){\vector(3,-2){23.5}}
 \thicklines
 \put(0,0){\vector(3,2){23.5}}
 \put(18,0){\vector(3,2){23.5}}
 \put(18,16){\vector(3,-2){23.5}}
 \put(24,16){\vector(3,-2){23.5}}
 \put(30,16){\vector(3,-2){23.5}}
 \end{picture}
 \end{center}

\medskip
The fact that the feasible functions one-to-one correspond to the vertices
of the crystal $K(\bfc)$ can be shown by two methods. A direct proof of
the assertion that $\Fscr$ along with the moves obeys axioms (A1)--(A5)
will be given in Section~\ref{sec:mod-cryst}. Another way consists in
showing a correspondence to GT-patterns and relies on
property~\refeq{cr-GT}. For $p,q\in\{1,\ldots,n\}$ with $p\le q$, let
$c[p:q]$ denote $c_p+\ldots+c_q$. As before, $c^\Sigma_j$ stands for
$c[1:n-j+1]$.

  \begin{prop} \label{pr:f-patt}
For $1\le j\le i\le n$, define
  \begin{equation} \label{eq:xij}
x_{i,j}:=\bar f_i(j)+c[1:i-j],
  \end{equation}
where $\bar f_i(j)$ denotes the sum of values of $f$ on the nodes in
$V_i(j)$. This gives a bijection between the set of feasible functions $f$
and the set of GT-patterns $X=(x_{i,j})$ of size $n$ bounded by
$\bfc^\Sigma$.
  \end{prop}

(Note that this leads to an alternative proof of~\refeq{cr-GT}, via the
crossing model.)

\medskip
  \begin{proof}
For a weight function $f$ satisfying~\refeq{feas}(i),(ii) (but not
necessarily~\refeq{feas}(iii)), define $X$ by~\refeq{xij}. Each multinode
$V_n(j)$ in the bottom level consists of the single node
$v=v_n^{n-j+1}(j)$, and we have $0\le f(v)\le c_{n-j+1}$ (since $v$ is in
$G^{n-j+1}$). Therefore, $x_{n,j}$ is between $c[1:n-j]$ and $c[1:n-j+1]$.

The inequality $x_{i,j}\ge x_{i+1,j+1}$ is provided by non-increasing $f$
along the edges from $V_i(j)$ to $V_{i+1}(j+1)$ and by the fact that the
term in~\refeq{xij} concerning $\bfc$ is the same for $(i,j)$ and for
$(i+1,j+1)$. The inequality $x_{i+1,j}\ge x_{i,j}$ follows from
non-increasing $f$ along the edges from $V_{i+1}(j)$ to $V_i(j)$ and from
the inequality $c_{i-j+1}\ge f(v_i^{i-j+1}(j))$. Thus, $X$ is a GT-pattern
bounded by $\bfc^\Sigma$.

Conversely, let $X$ be a GT-pattern bounded by $\bfc^\Sigma$. We construct
the desired $f$ step by step, starting from the bottom level. For each
node $v=v_n^{n-j+1}(j)$ (forming $V_n(j)$), we define
$f(v):=x_{n,j}-c[1:n-j]$. This value is nonnegative, and~\refeq{xij} holds
for $i=n$.

Now consider a multinode $V_i(j)$ with $i<n$, assuming that $f$ is already
determined for all levels $i'>i$ and satisfies~\refeq{feas}
and~\refeq{xij} for the nodes in these levels and the edges between them.
We show that $f$ can be properly extended to the nodes in $V_i(j)$ and
that such an extension is unique. Consider an {\em intermediate} node $v$
in $V_i(j)$ (existing when $i<n-1$). It has both SW- and SE-edges, say,
$(u,v),(v,w)$. The weights of $u$ and $w$ (already defined) satisfy
$f(u)\ge f(w)$ (since for the node $v'$ in the level $i+2$ such that
$(u,v')=\eSE(u)$ and $(v',w)=\eSW(w)$, we have $f(u)\ge f(v')\ge f(w)$).
The maximum possible weight of $v$ not violating~\refeq{feas}(i) is
$f(u)$, while the minimum possible weight is $f(w)$. In its turn, the {\em
first} node $v$ of $V_i(j)$ is connected with the level $i+1$ by the
unique edge $\eSE(v)$, say, $(v,w)$, and the maximum possible weight of
$v$ is $c_{i-j+1}$ (since $v$ belong to $G^{i-j+1}$), while the minimum
one is $f(w)$. And the {\em last} node $v$ of $V_i(j)$ is connected with
the level $i+1$ by the unique edge $\eSW(v)$, say, $(u,v)$, the maximum
possible weight of $v$ is $f(u)$, and the minimum one is zero.

Thus, the maximum assignment of weights for all nodes of $V_i(j)$) would
give $\bar f_i(j)=\bar f_{i+1}(j)+c_{i-j+1}$, implying $x_{i,j}\le \bar
f_i(j)+c[1:i-j]$, in view of  $x_{i,j}\le x_{i+1,j}=\bar
f_{i+1}(j)+c[1:i-j+1]$. And the minimum assignment would give $\bar
f_i(j)=\bar f_{i+1}(j+1)$, implying $x_{i,j}\ge \bar f_i(j)+c[1:i-j]$, in
view of $x_{i,j}\ge x_{i+1,j+1}= \bar f_{i+1}(j+1)+c[1:i-j]$. Therefore,
starting with the maximum assignment, scanning the nodes in $V_i(j)$
according to their ordering and decreasing their weights step by step, one
can always correct the weights so as to satisfy~\refeq{feas}(iii)
and~\refeq{xij}, while maintaining~\refeq{feas}(i),(ii). Moreover,
\refeq{feas}(iii) guarantees that the weights within $V_i(j)$ are
determined uniquely. Eventually, after handling leve1 1, we obtain the
desired function $f$ on $V(G)$.
  \end{proof}

\section{Moves in the model}
\label{sec:moves}

So far, we have dealt with the case of nonnegative {\em upper} bounds
(parameters) $c_1,\ldots,c_n$ and zero {\em lower} bounds, i.e., for any
feasible function $f$, the weight $f(v)$ of each node $v$ of a $k$-th base
subgraph lies between 0 and $c_k$. However, it is useful for us to
slightly extend the setting by admitting nonzero lower bounds (in
particular, for purposes of Subsection~\ref{sec:princ}.3 where the model
is extended to produce crystals with possible infinite monochromatic
paths).

Formally: for $\bfc,\bfd\in\Zset^n$ with $\bfc\ge \bfd$, we define a
feasible function to be an integer function $f$ on $V(G)$
satisfying~\refeq{feas}(i),(iii) and the relation
  \begin{equation}  \label{eq:bounds}
d_k\le f(v_i^k(j))\le c_k\qquad \mbox{for all $k,i,j$,}
  \end{equation}
instead of~\refeq{feas}(ii). The set of feasible functions for
$(\bfc,\bfd)$ is denoted by $\Fscr(\bfc,\bfd)$. Clearly the numerical part
of the model remains equivalent when for any $k$, we add a constant to
both $c_k$ and $d_k$ and accordingly add this constant to any weight
function for $G^k$. In particular, $\Fscr(\bfc,\bfd)$ is isomorphic to
$\Fscr(\bfc-\bfd,\bfzero)$, and when $\bfd=\bfzero$, $\Fscr(\bfc,\bfd)$
coincides with $\Fscr(\bfc)$ as above.

Now we start describing the desired transformations of functions in
$\Fscr(\bfc,\bfd)$, or moves (that will correspond to edges of the crystal
$K(\bfc-\bfd)$). Each transformation is performed only within one level
$i$, in which case it is called an $i$-{\em move}. We need some additional
definitions, notation and construction.

First of all, to simplify technical details, we extend each $G^k$ by
adding {\em extra} nodes and edges. More precisely, in the extended
digraph $\bar G^k$, the node set consists of elements $v^k_i(j)$ for
$(i,j)=(0,0)$ and for all $i,j$ such that $0\le i,j\le n+1$ and $j\le
i+1$, except for $(i,j)=(n+1,0)$. The edge set of $\bar G^k$ consists of
all possible pairs of the form $(v^k_i(j),v^k_{i-1}(j))$ or
$(v^k_i(j),v^k_{i+1}(j+1))$ (as before). An instance is illustrates in the
picture; here $n=4$, $k=2$, and the thick lines indicate the edges of the
original graph $G_4^2$.
  \begin{center}
  \unitlength=1mm
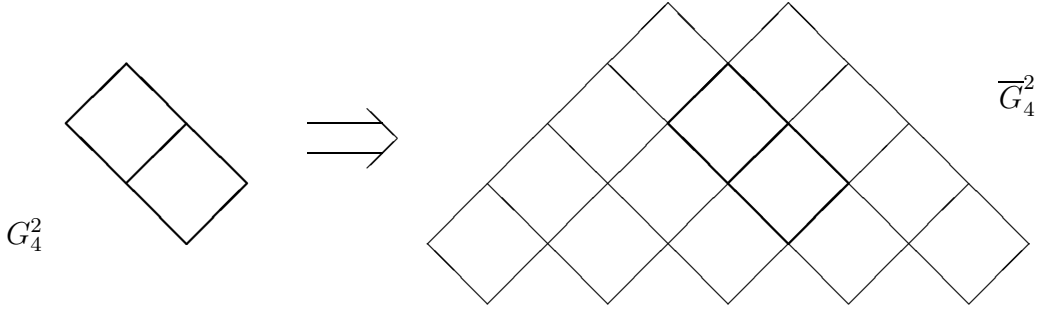
  \begin{picture}(140,40)
\put(56,8){\line(1,1){32}}
\put(64,0){\line(1,1){40}}
\put(80,0){\line(1,1){32}}
\put(96,0){\line(1,1){24}}
\put(112,0){\line(1,1){16}}
\put(128,0){\line(1,1){8}}
\put(64,0){\line(-1,1){8}}
\put(80,0){\line(-1,1){16}}
\put(96,0){\line(-1,1){24}}
\put(112,0){\line(-1,1){32}}
\put(128,0){\line(-1,1){40}}
\put(136,8){\line(-1,1){32}}
\put(40,20){\line(1,0){10}}
\put(40,24){\line(1,0){10}}
\put(52,22){\line(-1,1){4}}
\put(52,22){\line(-1,-1){4}}
  \thicklines
\put(8,24){\line(1,1){8}}
\put(16,16){\line(1,1){8}}
\put(24,8){\line(1,1){8}}
\put(24,8){\line(-1,1){16}}
\put(32,16){\line(-1,1){16}}
\put(88,24){\line(1,1){8}}
\put(96,16){\line(1,1){8}}
\put(104,8){\line(1,1){8}}
\put(104,8){\line(-1,1){16}}
\put(112,16){\line(-1,1){16}}
\put(0,8){$G^2_4$}
\put(132,26){$\bar G^2_4$}
  \end{picture}
 \end{center}
The disjoint union of these $\bar G^k$ gives the {\em extended supporting
graph} $\bar G$. It possesses the property that the original multinodes
become balanced, in the sense that for the set $J$ of index pairs $(i,j)$
satisfying $1\le j\le i\le n$, the extended multinodes $\bar V_i(j)$
contain the same number $n$ of nodes (these are
$v^1_i(j),\ldots,v^n_i(j)$). Also each node $v=v^k_i(j)$ of $\bar G$ with
$(i,j)\in J$ has exactly four incident edges, namely, all of $\eNW(v)$,
$\eSW(v)$, $\eNE(v)$, and $\eSE(v)$.

Each feasible function on $V(G)$ is extended to the extra nodes
$v=v^k_i(j)$ as follows:

(i) put $f(v):=c_k$ if there is a path from $v$ to $G^k$ (equivalently:
$j=0$ or $i-j>k-1$; one may say that $v$ lies on the left from $G^k$); and

(ii) put $f(v):=d_k$ otherwise (equivalently: $j>i$ or $j>n-k+1$, saying
that $v$ lies on the right from $G^k$).

One can see that such the extension maintains
conditions~\refeq{feas}(i),(ii),(iii) everywhere. Also
  \begin{numitem}
each edge of $\bar G$ with both ends not in $G$ is tight; and for any
$(i,j)\in J$, a node $v\in V_i(j)$ satisfies the switch condition in $G$
if and only if it does so in $\bar G$.
  \label{eq:switch_preserve}
  \end{numitem}

Given a feasible function $f$ on $V(\bar G)$, the move from $f$ in a level
$i\in \{1,\ldots,n\}$ changes $f$ within some multinode in this level. The
choice of this multinode depends on so-called residual slacks.

First, for a node $v=v^k_i(j)$, define
  $$
  \eps(v):=\partial f(\eNW(v)) \qquad \mbox{and} \qquad
  \delta(v):=\partial f(\eSE(v))
  $$
when the corresponding NW- or SE-edge exists in $\bar G$ (i.e., when
$i,j\ge 1$ in the former case and $i,j\le n$ in the latter case). We call
these the {\em upper slack} and the {\em lower slack} of $f$ at $v$,
respectively.

Next, define the upper slack $\eps_i(j)$  and the lower slack
$\delta_i(j)$ at a multinode $\bar V_i(j)$ as
  \begin{equation} \label{eq:slack}
   \eps_i(j):=\sum\nolimits_{k=1}^n \eps(v^k_i(j)) \qquad \mbox{and} \qquad
   \delta_i(j):=\sum\nolimits_{k=1}^n \delta(v^k_i(j))
  \end{equation}
(the former when $i,j\ge 1$, and the latter when $i,j\le n$). Note that
  \begin{equation} \label{eq:delta0}
  \eps_i(1)\ge \delta_i(0) \qquad \mbox{for $i=1,\ldots,n$}
   \end{equation}
(as $f(v^k_{i-1}(0))=f(v^k_i(0))=c_k$ and $f(v^k_i(1))\le f(v^k_{i+1}(1))$
for all $k$). Also $\eps_i(i+1)=0$.

Finally, we define the {\em residual upper slack} $\tilde\eps_i(j)$ and
the {\em residual lower slack} $\tilde\delta_i(j)$ by the following rule:
  \begin{numitem}
(a) for $0\le p<q\le i+1$, put $\pi(p,q):=\sum_{r=p+1}^q \eps_i(r)-
\sum_{r=p}^{q-1}\delta_i(r)$; (b) for $j=1,\ldots,i+1$, put
$\tilde\eps_i(j):=\max\{0,\min\{\pi(p,j)\colon 0\le p<j\}\}$; and (c) for
$j=0,\ldots,i$, put $\tilde\delta_i(j):=\max\{0,\min\{-\pi(j,q)\colon j<
q\le i+1\}\}$.
  \label{eq:reduce}
  \end{numitem}

  \noindent
{\bf Remark 2.} The residual slacks can also be computed via the following
recursive {\em cancelation process}. Initially, put $\tilde\eps_i:=\eps_i$
and $\tilde\delta_i:=\delta_i$. At each step, choose some pair $j'<j$ such
that $\tilde\eps_i(j)>0$, $\tilde\delta_i(j')>0$, and
$\tilde\eps_i(q)=\tilde\delta_i(q)=0$ for all $j'<q<j$. Subtract from each
of $\tilde\eps_i(j)$ and $\tilde\delta_i(j')$ their minimum. Repeat. Upon
termination of this process (when $j',j$ as above no longer exist), we
obtain $\tilde\eps$ and $\tilde\delta$ exactly as in~\refeq{reduce} for
all $j$. This observation will be used in Section~\ref{sec:mod-cryst}.

 \medskip
The residual slacks are integers and there exists $j\in\{1,\ldots,i\}$
such that
  \begin{equation}   \label{eq:activ}
\tilde\delta_i(0)=\ldots=\tilde\delta_i(j-1)=0\quad \mbox{and}\quad
\tilde\eps_i(j+1)=\ldots=\tilde\eps_i(i+1)=0.
  \end{equation}
(Indeed, suppose $\tilde\eps_i(j)>0$ and $\tilde\delta_i(j')>0$ for some
$j>j'$. Then, by~\refeq{reduce}(b),(c), $\pi(j',j)>0$ and $-\pi(j',j)>0$;
a contradiction.) Take the {\em minimum} $j$ satisfying~\refeq{activ} (if
there are many). If $\tilde\eps_i(j)>0$, then we say that $\bar V_i(j)$ is
the {\em active} multinode in the level $i$ (otherwise
$\tilde\eps_i(1)=\ldots=\tilde\eps_i(i)=0$ takes place).

The {\em moving operator} $\phi_i$ in level $i$ is applicable when the
active multinode $\bar V_i(j)$ does exist, and its action is simple: it
increases by one the value of $f$ on the switch-node in $\bar V_i(j)$,
preserving $f$ on all other nodes of $\bar G$.

To show that $\phi_i$ is well-defined, we will examine some rhombi of
$\bar G$, where by a (little) {\em rhombus} we mean a quadruple $\rho$ of
nodes of the form $v^k_i(j),v^k_{i-1}(j),v^k_i(j+1), v^k_{i+1}(j+1)$,
called the {\em left, upper, right, {\em and} lower} nodes of $\rho$,
respectively. The following simple observation is useful:
  \begin{numitem}
for a rhombus $\rho$, define $\partial f(\rho):=f(u')+f(w')-f(z')-f(v')$,
where $z',u',v',w'$ are the left, upper, right and lower nodes of $\rho$,
respectively; then $\partial f(\rho)=\eps(v')-\delta(z')$; in particular,
$\eps(v')-\delta(z')$ is nonnegative if $\eNE(z')$ or $\eSE(z')$ is
$f$-tight, and nonpositive if $\eNW(v')$ or $\eSW(v')$ is $f$-tight.
  \label{eq:simple}
  \end{numitem}

It follows that if $\bar V_i(j)$ is the active multinode in level $i$ and
$v$ is the switch-node in it, then $v$ belongs to $G$. Indeed, suppose $v$
occurs before the first node of $V_i(j)$. Then, in view
of~\refeq{feas}(iii) and~\refeq{switch_preserve}, the SW-edges of all
nodes of $\bar V_i(j)$ are tight, implying $\pi(j-1,j)\le 0$
(by~\refeq{simple}), contrary to $\tilde\eps_i(j)>0$. The fact that $v$
cannot occur after the last node of $V_i(j)$ is easy as well. So we can
speak of active multinodes within $G$ and of the switch-nodes there.

  \begin{prop} \label{pr:move}
The function $f':=\phi_i(f)$ is feasible.
  \end{prop}
  \begin{proof}
We have to check validity of (i) and (iii) in~\refeq{feas} for $f'$
(then~\refeq{bounds} for $f'$ will follow automatically). Below, when
speaking of switch-nodes or using expressions with
$\eps,\tilde\eps,\delta,\tilde\delta,\pi$, we always mean the
corresponding objects for $f$. Let $X:=\bar V_i(j)$ be the active
multinode for $f$ and $i$. Denote the nodes in $X$ by $v^1,\ldots,v^n$ (in
this order), and the switch-node by $v$.

Suppose $\partial f'(e)<0$ for some edge $e$. This is possible only if
$\partial f(e)=0$ and $e$ enters $v$, i.e., $e$ is $\eNW(v)$ or $\eSW(v)$.

\smallskip
(a) Let $e=\eSW(e)$. If $v\ne v^1$, then $\partial f(e)>0$ (otherwise the
switch-node in $X$ would occur before $v$). So $v=v^1$. Then the SW-edges
of all nodes in $X$ are $f$-tight, by~\refeq{feas}(iii). In view
of~\refeq{simple}, this implies $\pi(j-1,j)\le 0$, contrary to
$\tilde\eps_i(j)>0$.

\smallskip
(b) Now let $e=\eNW(v)$. The beginning node of $e$ belongs to the
multinode $V_{i-1}(j-1)$. Consider the rhombi $\rho^1,\ldots,\rho^n$
containing $v^1,\ldots,v^n$ as right nodes, respectively. Let
$z^k,u^k,w^k$ denote, respectively, the left, upper and lower nodes in
$\rho^k$. So $z^1,\ldots,z^n$ are the elements of $\bar V_{i}(j-1)$;~
$u^1,\ldots,u^n$ are the elements of $\bar V_{i-1}(j-1)$;~ and
$w^1,\ldots,w^n$ are the elements of $\bar V_{i+1}(j)$ (and the indices
grow according to the orderings in these multinodes). Let $v=v^p$, and let
$u^q$ be the switch-node in $\bar V_{i-1}(j-1)$. By~\refeq{feas}(iii), the
edges $(w^k,v^k)$ for $k=p+1,\ldots,n$ and the edges $(u^{k'},v^{k'})$ for
$k'=1,\ldots,q-1$ are tight for $f$. This gives
   $$
\eps(v^k)\le\delta(z^k)\qquad \mbox{for $k=1,\ldots,q-1$ and for
  $k=p+1,\ldots,n$}
   $$
(in view of~\refeq{simple}). Also the tightness of $e$ gives
$\eps(v^p)\le\delta(z^p)$. Suppose $q<p$. Then $u^p$ occurs in $\bar
V_{i-1}(j-1)$ after the switch-node $u^q$, and therefore, $(z^p,u^p)$ is
tight for $f$. We have $f(z^p)=f(u^p)=f(v^p)$, which implies the
$f$-tightness of all edges in $\rho^p$. Then $\partial f(\eSW(v))=0$,
contrary to shown in (a). Thus, $q\ge p$. This implies
$\eps(v^k)\le\delta(z^k)$ for all $k$, and therefore, $\eps_i(j)\le
\delta_i(j)$; a contradiction.

\smallskip
So, \refeq{feas}(i) for $f'$ is proven. Next, since $\partial
f'(e)\le\partial f(e)$ for all SW- and SE-edges $e$ of nodes in $\bar
V_{i-1}(j-1)$, \refeq{feas}(iii) is valid for $f'$ and this multinode.
Also~\refeq{feas}(iii) is, obviously, valid for $f'$ and $X$. It remains
to examine the multinode $Y:=V_{i-1}(j)$ since for the edge
$e=\eNE(v)=(v,u)$, which is the SW-edge for the node $u$ in $Y$, the value
$\partial f'(e)$ becomes greater than $\partial f(e)$. If $e$ is not
$f$-tight or if the {\em last} node $u'$ in $Y$ satisfying the switch
condition for $f$ does not occur before $u$, then~\refeq{feas}(iii)
follows automatically.

Suppose $\partial f(e)=0$ and $u'$ occurs before $u$. We show that this is
not the case by arguing in a way close to (b). For $k=1,\ldots,n$, let
$z^k,u^k,v^k,w^k$ denote, respectively, the left, upper, right and lower
nodes of the rhombus whose upper node (namely, $u^k$) is contained in $Y$.
Then the node $v$ (as before) is $z^p$, and $u'=u^q$ for some $p,q$ with
$q<p$. The fact that both $v,u'$ satisfy the switch condition for $f$ (in
their multinodes), together with $q<p$, implies that for each
$k=1,\ldots,n$, at least one of $\partial f(z^k,u^k)$ and
$\partial(z^k,w^k)$ is zero. This gives (cf.~\refeq{simple}):
    $$
  \delta(z^k)\le \eps(v^k)\qquad \mbox{for all $k$.}
    $$
Moreover, this inequality is strict for $k=q$. Indeed, we have $\partial
f(z^q,w^q)=0$ and $\partial f(u^q,v^q)>0$ (otherwise the node in $Y$ next
to $u$ would satisfy the switch condition for $f$ as well, but $u'$ is the
last of such nodes). So we obtain $\delta_i(j)<\eps_i(j+1)$. This implies
(in view of $\tilde\eps_i(j)>0$ and $\tilde\delta_i(j')=0$ for
$j'=0,\ldots,j-1$) that $\tilde\delta_i(j)=0$ and $\tilde\eps_i(j+1)>0$,
and therefore, the active multinode in level $i$ should occur after
$V_i(j)$; a contradiction.

This completes the proof of the proposition.
   \end{proof}

In conclusion of this section we discuss one more important aspect.

\medskip
\noindent {\bf Backward moves.} Besides the above description of partial
operators $\phi_i$ that increase functions in $\Fscr(\bfc,\bfd)$, we can
describe explicitly the corresponding decreasing operators, which make
{\em backward moves}. For $i=1,\ldots,n$, such an operator $\psi_i$ acts
on a feasible function $f$ as follows (as before, we prefer to deal with
extended functions on $V(\bar G)$). We take the {\em first} multinode
$V_i(j)$ (with $j$ minimum) in level $i$ for which $\tilde\delta_i(j)>0$;
the operator does not act when $\tilde\delta_i(j)=0$ for all $j$. In view
of~\refeq{delta0}, $1\le j\le n$.  In this multinode, called {\em active
in backward direction}, we take the {\em last} node $v$ possessing the
switch condition~\refeq{feas}(iii), called the {\em switch-node in
backward direction}. Then the action of $\psi_i$ consists in decreasing
the weight $f(v)$ by one, preserving the weights of all other nodes of
$G$.
  \begin{prop} \label{pr:backmove}
The function $f':=\psi_i(f)$ is feasible. Moreover,  $\phi_i$ is applicable to
$f'$, and $\phi_i(f')=f$.
  \end{prop}
  \begin{proof}
One can prove this by arguing in a similar spirit as in the proof of
Proposition~\ref{pr:move}. Instead, we can directly apply that proposition
to a certain reversed model. This is based on a simple observation, as
follows.

For a node $v\in V(\bar G)$, define $\mu(v):=\partial f(\eNE(v))$ and
$\nu(v):=\partial f(\eSW(v))$ (when such an NE- or SW-edge exists in $\bar
G$). The {\em alternative} upper and lower slacks at a multinode $V_i(j)$
are defined to be, respectively, the sum of numbers $\mu(v)$ and the sum
of numbers $\nu(v)$ for the nodes $v$ in this miltinode (the former is
defined for $j=0,\ldots,i$, and the latter for $j=1,\ldots,i+1$).
Compare~\refeq{slack}. Considering the little rhombus containing nodes
$u=v^k_i(j-1)$ and $v=v^k_i(j)$, we have $\nu(v)-\mu(u)=\eps(v)-\delta(u)$
(cf.~\refeq{simple}). This gives
  \begin{equation}  \label{eq:rhomb_rel}
\nu_i(j)-\mu_i(j-1)=\eps_i(j)-\delta_i(j-1).
  \end{equation}

The {\em reversed model} $\Mscr^r$ is obtained by reversing the edges of
$G$, by replacing the upper bound $\bfc$ by $-\bfd$, and by replacing the
lower bound $\bfd$ by $-\bfc$ (one may think that we now read the original
model from right to left). Accordingly, a feasible function $f$ in $M$ is
replaced by $f^r:=-f$. One can see that $f^r$ is feasible for $M^r$ and
that the {\em last} node satisfying the switch condition for $f$ in an
original multinode $V_i(j)$ turns into the switch-node for $f^r$ in the
corresponding multinode $V^r_i(j')$ in $M^r$. Also $\eps^r_i(j')=\mu_i(j)$
and $\delta^r_i(j')= \nu_i(j)$ (where $\eps^r,\delta^r$ stand for
$\eps,\delta$ in the reversed model). In view of~\refeq{rhomb_rel},
expressions in~\refeq{reduce} with $f^r$ in $M^r$ will give
$\tilde{\eps^r_i}(j')=\tilde\delta_i(j)$ and
$\tilde\delta^r_i(j')=\eps_i(j)$ for all $j$.

These observations enable us to conclude that the function $(f^r)'$
obtained by the forward move from $f^r$ in $M^r$ generates the function
$f'=\psi_i(f)$ in $M$. Therefore, $f'$ is feasible. To see the second part
of the proposition, let $v$ be the node of the active in backward
direction multinode $V_i(j)$ where $f$ decreases (by one) to produce $f'$.
The edge $\eSW(v)$ is non-tight for $f'$, which implies that $v$ is the
{\em unique} node in $V_i(j)$ satisfying the switch condition for $f'$,
and therefore, $v$ becomes the switch-node there. Also decreasing $f$ by
one at $v$ results in increasing $\eps(v)$, and one can see that the
residual slack $\tilde \eps_i(j)$ for $f'$ is greater by one than that for
$f$. This and~\refeq{activ} imply that $V_i(j)$ is just the active
multinode for $f'$ and $i$. Hence the forward move from $f'$ increases it
by one at $v$, and we obtain $\phi_i\psi_i(f)=f$, as required.
  \end{proof}

By this proposition, the operator $\psi_i$ is injective. The ``doubly
reversed'' model coincides with the original one, and therefore,
Proposition~\ref{pr:backmove} implies that $\psi_i\phi_i(f)=f$ for each
$f$ to which $\phi_i$ is applicable. So $\phi$ and $\psi$ are inverse to
each other and we may denote $\psi_i$ by $\phi_i^{-1}$.

\section{The relation of the model to RAN-crystals}
\label{sec:mod-cryst}

We have seen that the feasible functions in the model one-to-one
correspond to the vertices of a crystal, by using the GT-pattern model for
the latter, see Proposition~\ref{pr:f-patt}. In this section we directly
verify that the set $\Fscr$ of these functions and the set of (forward)
moves satisfies axioms (A1)--(A5), and therefore, they constitute an
RAN-crystal. One may assume that the lower bounds are zero, i.e.,
$\Fscr=\Fscr(\bfc)$ for $\bfc\in\Zset_+^n$. When the operator $\phi_i$ is
applicable to an $f\in \Fscr$, we say that $f$ and $f':=\phi_i(f)$ are
connected by the directed edge $(f,f')$ with the color $i$; the set of
these edges is denoted by $\Escr_i$. This produces the $n$-colored digraph
$\Kscr(\bfc)=(\Fscr,\Escr)$ in which $\Escr$ is partitioned into the color
classes $\Escr_1,\ldots,\Escr_n$. So we are going to show the following.

  \begin{theorem} \label{tm:mod-cryst}
$\Kscr(\bfc)$ is an RAN-crystal.
  \end{theorem}
  \begin{proof}
As before, it is more convenient to operate with the extended supporting
graph $\bar G$ and assume that the functions in $\Fscr$ are properly
extended to the nodes in $V(\bar G)-V(G)$.

Axiom (A1) immediately follows from properties of operators $\phi_i$ and
$\psi_i$. Next we observe the following. For $f\in\Fscr$ and a color $i$,
if $V_i(j)$ is the active multinode, then the action of $\phi_i$ decreases
$\tilde\eps_i(j)$ by 1, increases $\tilde\delta_i(j)$ by 1, and does not
change the residual slacks $\tilde\eps$ and $\tilde\delta$ for the other
multinodes in level $i$. This follows from~\refeq{activ} and the fact that
under increasing $f$ by 1 at the switch-node $v$ in $V_i(j)$, $\eps(v)$
decreases by 1 and $\delta(v)$ increases by 1. Similarly, if $V_i(j')$ is
the active multinode in backward direction, then $\psi_i$ decreases
$\tilde\delta_i(j')$ by 1, increases $\tilde\eps_i(j')$ by 1, and
preserves the residual slacks for the other multinodes in level $i$. This
implies
  \begin{equation}  \label{eq:th-in-level}
\ellout_i(f)=\sum\nolimits_{j=1}^i \tilde\eps_i(j) \qquad\mbox{and}
\qquad \ellin_i(f)=\sum\nolimits_{j=1}^i \tilde\delta_i(j),
  \end{equation}
regarding $f$ as a vertex of $\Kscr$.

If $i,i'$ are two colors with $|i-i'|\ge 2$, then any changes of $f$ in
the level $i$ do not affect the numbers $\eps(v)$ and $\delta(v)$ for
nodes $v$ in the level $i'$. So $\ellout_{i'}(f) =\ellout_{i'}(f')$ and
$\ellin_{i'}(f)=\ellin_{i'}(f')$ for $f'=\phi_i(f)$. This implies validity
of axiom (A5).

In order to verify axioms (A2),(A3) and (especially) (A4) for neighboring
colors, we need a more careful analysis of the behavior of residual
slacks. The following interpretation for the cancelation process (see
Remark~2 in Section~\ref{sec:moves}) is of help.

For $f\in\Fscr$ and a fixed level $i''$, we may think of $V(j)$ as a {\em
box} where $\eps(j)$ {\em white balls} and $\delta(j)$ {\em black balls}
are contained (we omit the subindex $i''$ hereinafter). Imagine that there
is a set $C$ of {\em couples}, each involving one black ball $b$ from a
box $V(j)$ and one white ball $w$ from a box $V(j')$ such that $j<j'$
(each ball occurs in at most one couple). We associate to a couple $(b,w)$
the integer {\em interval} $[j(b),j(w)]$, where
   $$
  \mbox{$j(q)$ denotes the number of the box containing a ball $q$}.
   $$
The set $\Iscr$ of these intervals (with possible multiplicities) is
required to form an {\em interval family}, which means that there are no
two intervals $[\alpha,\beta]$, $[\alpha',\beta']$ such that
$\alpha<\alpha'<\beta<\beta'$ (i.e., no crossing intervals). In
particular, the set of maximal intervals in $\Iscr$, not counting
multiplicities, forms a linear order in a natural way. Also it is required
that: (i) $C$ is maximal, in the sense that there are no uncoupled, or
{\em free}, a black ball $b$ and a white ball $w$ such that $j(b)<j(w)$;
and (ii) no free ball lies in the interior of an interval in $\Iscr$.

It is easy to realize that such a $C$ exists and unique, up to recombining
couples with equal intervals. We denote the set of free white (free black)
balls by $W$ (resp. $B$) and call $(C,W,B)$ the {\em arrangement}
for the given collection of black and white balls. Furthermore, for each $j$,
the number of free white balls (free black balls) in $V(j)$ is precisely
$\tilde\eps(j)$ (resp. $\tilde\delta(j)$).

Let $p$ denote the maximal number $j(w)$ among $w\in W$ (letting
$p=-\infty$ if $W=\emptyset$), and $q$ the minimal number $j(b)$
among $b\in B$ (letting $q=\infty$ if $B=\emptyset$). Then $p\le q$.
One can see that if some {\em black} ball $b$ is {\em removed},
then the arrangement changes as follows (we indicate only the changes
important for us).
  \begin{numitem}
If $b$ is free, it is simply deleted from $B$. And if $b$ is coupled and
occurs in a maximal interval $\sigma=[\alpha,\beta]$, then: (a) if
$\beta\le q$ then one of the previously coupled white balls $w$ with
$j(w)=\beta$ becomes free (and $\sigma$ is replaced by a maximal interval
$[\alpha,\beta']$ for some $j(b)<\beta'\le\beta$, unless $\sigma$ vanishes
at all); and (b) if $q\le \alpha$, then one free black ball $b'$ whose
number $j(b')$ is maximum provided that $j(b')\le \alpha$ becomes coupled
and generates the maximal interval $[j(b'),\beta]$.
   \label{eq:rem_b}
  \end{numitem}

On the other hand, when a new {\em white} ball $w$ is {\em added},
the changes are as follows.
   \begin{numitem}
In case $j(w)\le q$: (a) if $j(w)$ is in the interior of some maximal
interval $[\alpha,\beta]$, then $w$ becomes coupled and one
previously coupled white ball $w'$ with $j(w')=\beta$ becomes free;
(b) otherwise $w$ is simply added to $W$. And in case $j(w)>q$: (c)
$w$ becomes coupled and one free black ball $b$ with $j(b)$ maximum
provided that $j(b)<j(w)$ becomes coupled as well.
   \label{eq:add_w}
   \end{numitem}

Using this interpretation, we now check axioms (A2)--(A4) for neighboring
levels (viz. colors) $i$ and $i-1$ in the model. Here for $f\in\Fscr$ in
question, the number of the active multinode (the active multinode in
backward direction) in the level $i$ is denoted by $p=p(f)$ (resp.
$q=q(f)$), and $p'=p'(f)$ and $q'=q'(f)$ stand for the analogous numbers
in the level $i-1$ (as before, we use the sign $-\infty$ or $\infty$ if
such a multinode does not exist).

 \medskip
 \noindent {\bf Verification of (A2).}
When $\phi_i$ applies to $f$ (at $V_i(p)$), the value $\delta_{i-1}(p-1)$
decreases by 1. (Recall that for $v\in V_i(j)$ and $(u,v)=\eNW(v)$, $u$
belongs to $V_{i-1}(j-1)$.) In the above interpretation, this means that
one black ball is removed from the arrangement for the level $i-1$.
Then~\refeq{rem_b} implies that in case $p-1<q'$, the sum of values
$\tilde\eps_{i-1}(j)$ over $j$ (equal to $h_{i-1}(f)$) increases by 1,
while all $\tilde\delta_{i-1}(j)$ preserve. And if $p-1\ge q'$, then the
sum of values $\tilde\delta_{i-1}(j)$ (equal to $t_{i-1}(f)$) decreases by
1, while all $\tilde\eps_{i-1}(j)$ preserve. Also in the former case, we
obtain $p(f')\le p(f)$ and $q'(f')=q'(f)$, where $f':=\phi_i(f)$, and
therefore, the next application of $\phi_i$ will fall in the former case
as well (further increasing $h_{i-1}$). Next, when $\phi_{i-1}$ applies to
$f$, we observe from~\refeq{add_w} that: in case $p'\le q-1$, the sum of
$\tilde\eps_i(j)$ increases by one, while all $\tilde\delta_i(j)$
preserve, and in case $p'\ge q$, the sum of $\tilde\delta_i(j)$ decreases
by one, while all $\tilde\eps_i(j)$ preserve. Also in the former case,
$p'(f')\le p'(f)$ and $q(f')=q(f)$, where $f':=\phi_{i-1}(f)$, so the next
application of $\phi_{i-1}$ increases $h_i$ as well.

\medskip
\noindent {\bf Verification of (A3).} This is also easy. Let
$f':=\phi_i(f)$ and $f'':=\phi_{i-1}(f)$. Suppose $(f,f')$ has label 0.
Then $p-1\ge q'$ and $p'(f')=p'(f)$ (see the previous verification).
Moreover, the switch-node $u$ in $V_{i-1}(p')$ for $f$ remains the
switch-node for $f'$. (Indeed, since $p'\le p-1$, the slacks of the
SW-edges of all nodes in $V_{i-1}(p')$ preserve, and the slacks of their
SE-edges do not increase.) In its turn, $p'\le p-1\le q-1$ implies that
$(f,f'')$ has label 1, as required in the axiom. Also neither the active
multinode in the level $i$ nor the switch-node $v$ in it can change when
$\phi_{i-1}$ applies to $f$. Thus, both $\phi_{i-1}\phi_i$ and
$\phi_i\phi_{i-1}$ increase the original function $f$ by 1 on the same
elements $u,v$. A verification of the relation
$\phi_{i-1}\phi_i(f)=\phi_i\phi_{i-1}(f)$ in the case when $(f,f'')$ has
label 0 is similar.

\medskip
\noindent {\bf Verification of (A4).} This is somewhat more involved.
Assuming that both $\phi_i$ and $\phi_{i-1}$ are applicable to a
feasible function $f$, define $f_1:=\phi_i(f)$ and
$g_1:=\phi_{i-1}(f)$, and let both $(f,f_1)$ and $(f,g_1)$ have label
1. Then $p-1<q'$ and $p'+1\le q$ (where $p=p(f)$, and similarly for
$q,p',q'$).

Since $\ell(f,f_1)=1$, we have $h_{i-1}(f_1)=h_{i-1}(f)+1\ge 2$. Therefore,
we can define $f_2:=\phi_{i-1}(f_1)$ and $f_3:=\phi_{i-1}(f_2)$.
Similarly, we can define $g_2:=\phi_i(g_1)$ and $g_3:=\phi_i(g_2)$. Our
aim is to show that $\phi_i$ is applicable to $f_3$, that $\phi_{i-1}$ is
applicable to $g_3$, and that $\phi_i(f_3)=\phi_{i-1}(g_3)$. Two cases
are possible: $p'\le p-1$ and $p'\ge p$.

 \medskip
 \noindent
{\bf Case $p'\le p-1$.}~ For $k=1,2,3$, we denote
$p(f_k),q(f_k),p'(f_k),q'(f_k)$ by $p_k,q_k,p'_k,q'_k$, respectively;
similar numbers for $g_k$ are denoted by $\bar p_k,\bar q_k,\bar p'_k,\bar
q'_k$. We use the above interpretation and associate to each current
function $f'$ the corresponding arrangement $(C=C(f'),W=W(f'),B=B(f'))$ in
the level $i$ and the corresponding arrangement
$(C'=C'(f'),W'=W'(f'),B'=B'(f'))$ in the level $i-1$.

Since $\tilde\eps_i(p)>0$, there is a white ball $w\in W(f)$ with
$j(w)=p$. In view of $p-1<q'$, $w$ corresponds to a coupled black ball
$b'$ with $j(b')=p-1$ in the level $i-1$; let $[\alpha',\beta']$ be the
maximal interval for $C'(f)$ that contains $b'$. Then $p-1<\beta' \le q'$.
We also define the number $\beta$ as follows: if the point $p'+1$ lies in
the interior of some maximal interval $[\tilde\alpha,\tilde\beta]$ for
$C(f)$, put $\beta:=\tilde\beta$; otherwise put $\beta:=p'+1$. (The
meaning of $\beta$ is: in view of $p'+1\le q$, if a new white ball $\hat
w$ with $j(\hat w)=p'+1$ is added in the level $i$, then the arrangement
in this level changes so that there appears a free ball $w'$ with
$j(w')=\beta$; see~\refeq{add_w}.) Appealing to the interpretation, we can
precisely characterize the changes of
$\tilde\eps_i,\tilde\delta_i,\tilde\eps_{i-1}, \tilde\delta_{i-1}$ when
the above-mentioned transformations of our functions are carried out.

\smallskip
(i) The transformation $f\to f_1$ decreases $\tilde\eps_i(p)$ by 1
and increases $\tilde\delta_i(p)$ by 1. Also
$\tilde\eps_{i-1}(\beta')$ becomes equal to 1; cf.~\refeq{rem_b}(a).

\smallskip
In particular, $p'_1=\beta'$, i.e., $V_{i-1}(\beta')$ becomes the active
multinode in the level $i-1$.

\smallskip
(ii) The transformation $f_1\to f_2$ reduces
$\tilde\eps_{i-1}(\beta')$ to 0 and increases
$\tilde\delta_{i-1}(\beta')$ by 1. Also $\tilde\delta_i(r)$
decreases by 1 for some $r\ge p=q_1$; cf.~\refeq{add_w}(c).

\smallskip
This gives $p'_2=p'$ and $q_2\ge p$ and preserves all intervals for $C$
that lie before $p$.

\smallskip
(iii) The transformation $f_2\to f_3$ decreases $\tilde\eps_{i-1}(p')$ by
1 and changes $\tilde\delta_{i-1}(p')$ from 0 to 1. Also
$\tilde\eps_i(\beta)$ increases by 1; cf.~\refeq{add_w}(a),(b).

\smallskip
The latter property implies $p'+1\le\beta\le p_3\le p$. Then $\phi_i$ is
applicable to $f_3$; define $f_4:=\phi_i(f_3)$. (Furthermore, one can see
that $V_i(p_3)$ is the active multinode in the level $i$ for the function
$\phi_i\phi_{i-1}(f)$ as well.)

Thus, the combined transformation $\phi_i\phi_{i-1}\phi_{i-1}\phi_i$
consecutively increases $f$ by 1 in the switch-nodes $v_0,v_1,v_2,v_3$ of
$V_i(p),V_{i-1}(\beta'),V_{i-1}(p'),V_i(p_3)$, respectively, where each
switch-node is defined for the current function at the moment of the
corresponding transformation. (Note that $p'$ and $\beta'$ are different,
while $p$ and $p_3$ may coincide.)

Next we examine the other chain of transformations.

\smallskip
(iv) The transformation $f\to g_1$ decreases $\tilde\eps_{i-1}(p')$ by 1
and increases $\tilde\delta_{i-1}(p')$ by 1. Also $\tilde\eps_i(\beta)$
increases by 1.

\smallskip
From~\refeq{add_w}(a),(b) it follows that $\beta\le p$, implying
$\bar p_1=p$.

\smallskip
(v) The transformation $g_1\to g_2$ decreases $\tilde\eps_i(p)$ by 1 and
increases $\tilde\delta_i(p)$ by 1. Also $\tilde\delta_{i-1}(p')$ reduces to 0.

\smallskip
Moreover,~\refeq{rem_b}(b) implies the following important property
($\ast$): $[p',\beta']$ becomes a maximal interval in the new arrangement
in the level $i-1$. Also (as mentioned after~(iii)) $\bar p_2$ coincides
with $p_3$.

\smallskip
(vi) The transformation $g_2\to g_3$ decreases $\tilde\eps_i(\bar p_2)$ by
1 and increases $\tilde\delta_i(\bar p_2)$ by 1. Also, in view of
$p'+1\le\bar p_2\le p$, the interval $[p',\beta']$ in the level $i-1$
(see~($\ast$) above) is destroyed and $\tilde\eps_{i-1}(\beta')$ becomes
equal to 1; cf.~\refeq{rem_b}(a).

\smallskip
So $\bar p'_3=\beta'$ and we can apply $\phi_{i-1}$ to $g_3$; let
$g_4:=\phi_{i-1}(g_3)$. We assert that $g_4=f_4$.

To see this, notice that the combined transformation
$\phi_{i-1}\phi_i\phi_i\phi_{i-1}$ increases the initial $f$ within the
same multinodes as those in the transformation
$\phi_i\phi_{i-1}\phi_{i-1}\phi_i$, namely, $V_{i-1}(p'),V_i(p),V_i(\bar
p_2=p_3),V_{i-1}(\beta')$ (but now the order is different). Let $\bar
v_0,\bar v_1,\bar v_2,\bar v_3$ be the switch-nodes in these multinodes,
respectively (each being taken at the moment of the corresponding
transformation). Since no change in the level $i-1$ affects the slacks of
SW- and SE-edges in the level $i$, we have $\bar v_1=v_0$ and $\bar
v_2=v_3$. Also $p'+1\le\beta,p$ implies that the transformations in the
level $i$ do not decrease the slacks of the SE-edges of nodes in
$V_{i-1}(p')$ and do not change the slacks of their SW-edges, whence $\bar
v_0=v_2$.

It remains to check that $\bar v_3=v_1$. Let $u$ be the switch-node in
$V_{i-1}(\beta')=:X$ for the initial function $f$. We have $p\le\beta'$.
Therefore, the increase at $v_0=\bar v_1$ can change the switch-node in
$X$ only if $p=\beta'$ and if the end $u'$ of the edge $\eNE(v_0)$ is
situated after $u$ in the ordering on $X$. If this is the case, then under
each of the transformations $f\to f_1$ and $g_1\to g_2$ (concerning
$V_i(p)$) the switch-node $u$ in $X$ is replaced by $u'$. Besides these,
there is only one transformation in the level $i$ that preceedes the
transformation within $X$, namely, $g_2\to g_3$. We know that $\bar p_2\le
p$ and that if $\bar p_2=p$ then $\bar v_2$ coincides with or preceedes
$v_0$ (taking into account that the transformation $g_1\to g_2$ concerning
$V_i(p)$ was applied earlier). This easily implies that $g_2\to g_3$ can
never change the switch-node in $X$. Thus, $\bar v_3=v_1$.

\smallskip
The case $p'\ge p$ is examined in a similar fashion, and we leave it to
the reader.

Finally, due to Proposition~\ref{pr:A4}, verifying the second part
of axiom (A4) (concerning the operators $\phi_i^{-1}$ and $\phi_{i-1}^{-1}$)
is not necessary.

This completes the proof of Theorem~\ref{tm:mod-cryst}.
  \end{proof}

\noindent {\bf Remark 3.} In light of the second claim in
Proposition~\ref{pr:A4}, instead of the tiresome verification of axiom
(A4) in the above proof, one may attempt to show that a maximal connected
subgraph with colors $i$ and $i-1$ in $\Kscr$ has only one zero-indegree
vertex. However, no direct method to show this is known to us.

\medskip
Clearly the source of the crystal $\Kscr(\bfc)$ is the identically
zero function $f_0$ on $V(G)$, and the sink is the function $f_\bfc$
taking the constant value $c_k$ within each subgraph $G^k$,
$k=1,\ldots,n$. In particular, this implies that
  \begin{numitem}
the distance (viz. the number of edges of a path) from the source to the
sink, or the {\em length} of $\Kscr(\bfc)$, is equal to
$\sum_{k=1}^n c_k|V(G^k)|$, or $\sum_{k=1}^n c_k k(n-k+1)$.
  \label{eq:distance}
  \end{numitem}

Also one can see that for the source function $f_0$ and a level $i$, one has
$\tilde\eps_i(1)=c_i$ and $\tilde\eps_i(j)=0$ for $j=2,\ldots,i$ (moreover:
starting from $f_0$, each application of $\phi_i$ increases the weight of
$v_i^i(1)$ by 1 until the weight becomes $c_i$). So $h_i(f_0)=c_i$ for each
color $i$. This means that $\Kscr(\bfc)$ is the crystal $K(\bfc)$, and now the
result of Stembridge~\cite{ste-03} that there exists exactly one RAN-crystal
with source having a prescribed $n$-tuple $\bfc$ of parameters
(see~\refeq{sourceAn}) and Corollary~\ref{cor:RANsource} enable us to conclude
with the following

  \begin{theorem} \label{tm:equiv}
The crossing model $\Mscr_n$ generates precisely the set of
regular $A_n$-crystals.
  \end{theorem}

\section{Principal lattice, principal subcrystals, and skeleton}
  \label{sec:princ}

In this section we apply the crossing model to establish certain
structural properties of RAN-crystals. We consider the initial setting for
the crossing model, i.e., when the upper bounds are nonnegative integers
and the lower bounds are zeros. So we deal with a parameter tuple
$\bfc=(c_1,\ldots,c_n)\in\Zset_+^n$ and the set $\Fscr(\bfc)$ of feasible
functions in the model. As before, $G=(V(G),E(G))$ is the supporting
graph, and $G^k=(V(G^k),E(G^k))$ is $k$-th base subgraph (component) in
$G$. The pair $(\Fscr(\bfc),\Escr(\bfc))$ is isomorphic to the crystal
$K=K(\bfc)=(V,E)$. Recall that $F_i$ denotes $i$-th partial operator on
$V$ (corresponding to the partial operator $\phi_i$ on $\Fscr(\bfc)$), and
$K^k= K^k_n(c_k)$ denotes $k$-th base crystal. We will also use the
following additional notation:

$v(f)$ denotes the vertex of $K$ corresponding to a feasible function $f$;

$f^1\sqcup f^2\sqcup\ldots\sqcup f^n$, where $f^k:V(G^k)\to\Zset$
($i=1,\ldots,n$), denotes the function on $V(G)$ coinciding with $f^k$ within
each $G^k$;

$v^k(f^k)$ denotes the vertex of $K^k$ corresponding to a feasible function
$f^k$ on $V(G^k)$;

$C^k a$ denotes the function on $V(G^k)$ taking a constant value
$a\in\Zset$.

 \bigskip
 \noindent
 {\bf\large \arabic{section}.1. Principal lattice and principal
 subcrystals.}~
Among the variety of feasible functions, certain functions are of most interest
to us. These are functions $f$ of the form $C^1a_1\sqcup\ldots\sqcup C^na_n$,
where each $a_k$ is an integer satisfying $0\le a_k\le c_k$. Such an $f$ is
feasible (since all edges of $G$ are $f$-tight); we call it a {\em principal
function} and denote by $f[\bfa]$, where $\bfa=(a_1,\ldots,a_n)$. The
corresponding vertex $v(f)$ is called a {\em principal vertex} of the crystal
and denoted by $v[\bfa]$. In particular, the source and sink of $K$ are the
principal vertices $v[\bfzero]$ and $v[\bfc]$, respectively. So there are
$(c_1+1)\times\ldots\times(c_n+1)$ principal vertices; their set is denoted by
$\Pi=\Pi(\bfc)$ and called the {\em principal lattice} in $K$.

The principal lattice possesses a number of nice properties, described
throughout this and next sections. One of them is that the intervals
between pairs of principal vertices are RAN-crystals as well, where for
vertices $u,v$ in an (acyclic) digraph, the {\em interval} from $u$ to $v$
is the subgraph $\Int(u,v)$ formed by the vertices and edges lying on
paths from $u$ to $v$.

To show this (and also for purposes of Subsection~\arabic{section}.3), we
first consider the crossing model with tuples $\bfc',\bfd'\in\Zset^n$ of
upper and lower bounds, $\bfc'\ge\bfd'$. This gives the crystal
$K(\bfc'-\bfd')$, also denoted as $K(\bfc',\bfd')$. Let
$\bfc'',\bfd''\in\Zset^n$ be such that $\bfc''\ge\bfc'$ and
$\bfd''\le\bfd'$. Clearly
  \begin{numitem}
any feasible function $f$ for $(\bfc',\bfd')$ is feasible for
$(\bfc'',\bfd'')$ as well.
   \label{eq:cdc'd'}
  \end{numitem}

This gives an injective map $\gamma$ from the vertex set of
$K'=K(\bfc',\bfd')$ to the vertex set of the crystal
$K''=K(\bfc'',\bfd'')$. Comparing the residual slacks $\tilde\eps'_i(j)$
and $\tilde\delta'_i(j)$ for the function $f$ in the model with the bounds
$\bfc',\bfd'$ and the residual slacks $\tilde\eps''_i(j)$ and
$\tilde\delta''_i(j)$ for $f$ in the model with the bounds
$\bfc'',\bfd''$, one can see that
  \begin{eqnarray}
 \tilde\eps''_i(1)=\tilde\eps'_i(1)+c''_i-c'_i\quad\mbox{and}\quad
\tilde\eps''_i(j)=\tilde\eps'_i(j)\quad &\mbox{for}& \;\; j=2,\ldots,i;
                  \label{eq:eps-eps'}\\
  \tilde\delta''_i(i)=\tilde\delta'_i(i)+d'_i-d''_i\quad\mbox{and}\quad
\tilde\delta''_i(j')=\tilde\delta'_i(j')\quad &\mbox{for}& \;\;
                         j'=1,\ldots,i-1. \nonumber
  \end{eqnarray}
Moreover, for each multinode, the switch-nodes concerning $f$ in both
models are the same, and similarly for the switch-nodes in backward
direction. Also the situation when an active multinode $V_i(j)$ for
$f,\bfc'',\bfd''$ is not active for $f,\bfc',\bfd'$ can arise only if:
$j=1$, the switch-node in $V_i(j)$ is $v_i^i(1)$, and $f(v_i^i(1))=c_i$;
and symmetrically for the active multinodes in backward direction. These
observations show that $\gamma$ is extendable to the edges of $K'$, and
moreover,
  \begin{numitem}
the image of $K'$ by $\gamma$ is a subcrystal of $K''$ isomorphic to $K'$,
and any path in $K''$ connecting vertices of $\gamma(K')$ is entirely
contained in $\gamma(K')$. Therefore, $\gamma(K')$ is the interval
$\Int(\gamma(s_{K'}),\gamma(\bar s_{K'}))$ of $K''$, where $s_{K'}$ and
$\bar s_{K'}$ are the source and sink of $K'$, respectively.
   \label{eq:K-in-K'}
   \end{numitem}
Note that $\gamma(s_{K'})$ and $\gamma(\bar s_{K'})$ are the principal
vertices $v[\bfd'-\bfd'']$ and $v[\bfc'-\bfd'']$ in $K(\bfc''-\bfd'')$,
respectively. So we obtain the following
  \begin{prop} \label{pr:princ_int}
For $\bfa,\bfb,\bfc\in\Zset_+^n$ with $\bfa\le\bfb\le\bfc$, the interval
of $K(\bfc)$ between the principal vertices $v[\bfa]$ and $v[\bfb]$ is
isomorphic to the RAN-crystal $K(\bfb-\bfa)$.
  \end{prop}

 \bigskip
 \noindent
 {\bf\large \arabic{section}.2. Skeleton.}~
This is a certain part of a RAN-crystal $K=K(\bfc)$ related to so-called {\em
1-relaxations} of principal functions. We use notation $\bfa^{(-k)}$ for an
$(n-1)$-tuple of integers $a_i$ where the index $i$ ranges
$1,\ldots,k-1,k+1,\ldots,n$. For $\bfa^{(-k)}$ satisfying
$\bfzero^{(-k)}\le\bfa^{(-k)}\le\bfc^{(-k)}$, define $\Fscr[\bfa^{(-k)}]$ to be
the set of all feasible functions $f=f^1\sqcup\ldots\sqcup f^n$ on $V(G)$ such
that $f^i=C^ia_i$ for each $i\ne k$. In other words, the non-fixed part $f^k$
of $f$ is any feasible function for $G^k$. (The latter is an arbitrary
nonnegative integer function $g$ on $V(G^k)$ bounded by $c_k$ and satisfying
the monotonicity condition $g(u)\ge g(v)$ for each edge $(u,v)\in E(G^k)$.
Since the switch condition becomes redundant for $G^k$ taken separately, just
all these functions $g$ generate the vertices $v$ of $K^k$: $v=v^k(g)$.)

Let $K[\bfa^{(-k)}]$ denote the subgraph of $K$ induced by the set of
vertices $v(f)$ for all $f\in\Fscr[\bfa^{(-k)}]$. For any
$f\in\Fscr[\bfa^{(-k)}]$, all edges in the subgraphs $G^i$ with $i\ne k$
are $f$-tight. Also each multinode $X$ of $G$ contains at most one node of
$G^k$. These facts imply that the moves from $f$ do not depend on the
entries of $\bfa^{(-k)}$, unless $f$ is transformed within a leftmost
multinode $V_i(1)$ not intersecting $G^k$ (i.e., with $i\ne k$). This
leads to the following property.

  \begin{prop} \label{pr:isom-subgr}
For any $\bfa^{(-k)}\le\bfc^{(-k)}$, the subgraph $K[\bfa^{(-k)}]$ of
$K(\bfc)$ is isomorphic to the base crystal $K^k(c_k)$.
  \end{prop}

The union $\Cscr$ of these subgraphs $K[\bfa^{(-k)}]$ over all $k$ and all
$\bfa^{(-k)}\le\bfc^{(-k)}$ constitutes the object that we call the {\em
skeleton} of $K$. Each $K[\bfa^{(-k)}]$ contains $c_k+1$ principal
vertices $v[\bfa']$ of $K(\bfc)$; here $a'_i=a_i$ for $i\ne k$ and $a'_k$
runs $0,\ldots,c_k$. The corresponding set of $c_k+1$ vertices in the base
crystal $K^k_n(c_k)$ is referred as its {\em axis} and denoted by
$S^k=S^k(c_k)$. (In case $n=2$, \cite{DKK-07} uses the name ``diagonal''
rather than ``axis''.)

The proposition below asserts that the skeleton of $K(\bfc)$ is obtained
from the base crystals $K^k(c_k)$ by use of a construction which is a
natural generalization of the diagonal-product construction for
RA2-crystals (see Theorem~\ref{tm:main2}) to the case of $n$ colors.

Again (like for $n=2$) we can describe such a construction for arbitrary
graphs $H^1,\ldots,H^n$ with distinguished vertex subsets $S^1,\ldots,S^n$
(respectively). Let $\Vscr$ be the collection of all $n$-element sets
containing exactly one vertex from each $S^i$. For $k=1,\ldots,n$, let
$\Vscr^{(-k)}$ be the collection of all $(n-1)$-element sets containing
exactly one vertex from each $S^i$ with $i\ne k$. For each $k$, take
$|\Vscr^{(-k)}|$ copies of $H^k$, each being indexed as $H^k_q$ for
$q\in\Vscr^{(-k)}$. We glue these copies together by identifying, for each
$q=\{v_1,\ldots,v_n\}\in\Vscr$ (where $v_k\in S^k$), the copies of
vertices $v_k$ in $H^k_{q\setminus\{v_k\}}$, $k=1,\ldots,n$, into one
vertex. The resulting graph is denoted as
$(H^1,S^1)$\sprod$\cdots\,$\sprod$(H^n,S^n)$.

In our case we take as $H^k$ the base crystal $K^k_n(c_k)$, and as the
distinguished subset $S^k$ the axis in it. The graph
$(H^1,S^1)$\sprod$\cdots\,$\sprod$(H^n,S^n)$ is called the {\em
axis-product} and denoted by $\Ascr(\bfc)$ (this is an $n$-colored digraph
where the edge colors are inherited from the base crystals). The principal
vertices in $\Ascr(\bfc)$ are defined to be those obtained by gluing
together vertices from the axes of graphs $K^k$. So the principal vertices
of $\Ascr(\bfc)$ one-to-one correspond to the principal functions in the
model, or to the $n$-tuples $\bfa\in\Zset^n_+$ with $\bfa\le\bfc$.

Summing up the above explanations, we have the following

  \begin{prop} \label{pr:karkas}
$K(\bfc)$ contains an induced subgraph $K'$ isomorphic to
$\Ascr(\bfc)=(H^1,S^1)$\sprod$\cdots\,$\sprod$(H^n,S^n)$ (respecting edge
colors). Moreover, $K'$ is determined uniquely and its vertices correspond to
the feasible functions $f^1\sqcup\ldots\sqcup f^k$ for $(G,\bfc)$ such that
each $f^i$ is a constant function on $V(G^i)$, except possibly for one function
$f^k$, which is an arbitrary feasible function for $(G^k,c_k)$.
   \end{prop}

Here the uniqueness can be shown as follows. The length of a path in $K$
from the source $v[\bfzero]$ to the sink $v[\bfc]$ is equal to
$\sum_{k=1}^n c_k|V(G^k)|)$ (see~\refeq{distance}). The length of a path
from the source to the sink in $\Ascr(\bfc)$ is the same. Therefore (since
$K$ is graded), the source of $K'$ must be at $v[\bfzero]$ and the sink of
$K'$ must be at $v[\bfc]$. Now it is easy to realize that $K'$ is
reconstructed in $K$ in a unique way.

Next, for two principal vertices $v[\bfa]$ and $v[\bfb]$, let us say that
the latter is the $k$-th {\em immediate successor} of the former if
$b_k=a_k+1$ and $a_i=b_i$ for all $i\ne k$. One can see that any possible
transformation of the function $f[\bfa]$ into $f[\bfb]$ (by use of forward
moves in the model) consists of a sequence of $|V(G^k)|$ moves, and the
corresponding sequence of nodes where the current function changes forms a
{\em linear order} on $V(G^k)$ agreeable with the poset structure of
$G^k$. In other words, this is an ordering $(u_1,\ldots,u_d)$ of the nodes
of $G^k$ such that for each $p=1,\ldots,d$, the set
$U_p=\{u_1,\ldots,u_p\}$ is an ideal in $G^k$ (i.e., no edge goes to $U_p$
from the complement). Each $U_p$ determines the function $g_p$ on $V(G^k)$
taking the value $a_k+1$ within $U_p$, and $a_k$ on the rest. Let
$\ell(p)$ denote the level number of $u_p$ in $G$, and let $f_p$ denote
the function on $V(G)$ formed from $f[\bfa]$ by replacing $C^ka_k$ on
$V(G^k)$ by $g_p$. One can check that $f_p$ coincides with the function
obtained from $f_{p-1}$ by the move in the level $\ell(p)$ (which just
increases the weight of $u_p$ by one). Thus, we have the following

  \begin{prop}  \label{pr:fundam}
For $k=1,\ldots,n$ and a principal vertex $v$ of $K$, if the $k$-th
immediate successor $w$ of $v$ exists, then each paths from $v$ to $w$ in
$K$ one-to-one corresponds to a linear order $(u_1,\ldots,u_d)$ for $G^k$
(where $d=|V(G^k)|$). Under this correspondence, the node $w$ can be
expressed as $F_{\ell(d)}\cdots F_{\ell(1)}v$, where $\ell(p)$ is the
level number of $u_p$.
  \end{prop}

For $k=1,\ldots,n$, the set of strings $\ell(d)\cdots \ell(1)$ as in this
proposition is denoted by $FS_n(k)$; this is invariant for all principal
vertices having the $k$-th immediate successor. We refer to any of such
strings as a {\em fundamental} one. As a special case, $FS_n(k)$ contains
the fundamental string
  \begin{numitem}
$S_{n,k}=w_{n,k,n-k+1}\cdots w_{n,k,2} w_{n,k,1}$, where the substring
$w_{n,k,i}$ is of the form $(i)(i+1)\cdots (i+k-1)$.
  \label{eq:fstring}
  \end{numitem}
(This corresponds to a route in $G^k$ (according to which the weights of nodes
are consecutively increased by 1) consisting of $n-k+1$ paths, as follows. We
starts from the source $\lvert^k$ and go as long as possible in the
NE~direction, up to the topmost node $\topvert^k$ (obtaining the string
$w_{n,k,1}$ of levels). Then we begin at the next node on the SW-side of $G^k$,
namely, $v^k_{k+1}(2)$, and again go in the NE~direction (yielding
$w_{n,k,2}$), and so on. At the final stage, we begin at the last node on the
SW-side, namely, $\botvert^k$, and go up to the sink $\rvert^k$ (yielding
$w_{n,k,n-k+1}$).)

\medskip
\noindent {\bf Example}. Let $n=3$. Since the graph $G^1$ forms a path, there
is only one fundamental string for $k=1$, namely, $321$. Similarly, $FS_3(3)$
consists of a unique string, namely, $123$. The set $FS_3(2)$ for the graph
(rhombus) $G^2$ consists of two strings: $2312$ and $2132$.

 \bigskip
 \noindent
 {\bf\large \arabic{section}.3. Infinite crystals.}~
So far, we have dealt with $n$-colored crystals having a finite set of
vertices, or finite crystals. However, by use of the crossing model one
can generate infinite analogs of RAN-crystals (arising when we admit
infinite monochromatic paths). Some applications of ``crystals'' of this
sort are indicated in~\cite{kas-02} in connection with modified quantized
enveloping algebras. Infinite analogs of RA2-crystals are discussed
in~\cite[Sec.~6]{DKK-07}.)

To obtain infinite RAN-crystals, we use the crossing model with
double-sided bounds and consider an upper bound
$\bfc\in(\Zset\cup\{\infty\})^n$ and a lower bound
$\bfd\in(\Zset\cup\{-\infty\})^n$ with $\bfc\ge\bfd$. More strictly: for a
variable $M\in\Zset_+$ and each color $i$, define $c^M_i$ to be $c_i$ if
$c_i<\infty$, and $\max\{M,d_i\}$ otherwise, and define $d^M_i$ to be
$d_i$ if $d_i>-\infty$, and $\min\{-M,c_i\}$ otherwise. When $M$ grows,
there appears a sequence of finite crystals $K(\bfc^M,\bfd^M)$, each
containing the previous crystal $K(\bfc^{M-1},\bfd^{M-1})$ as a principal
interval, by~\refeq{K-in-K'}. At infinity we obtain the desired
(well-defined) ``infinite crystal'' $K(\bfc,\bfd)$ (when $\bfc$ or/and
$\bfd$ is not finite).

Some trivial consequences of this construction are as follows. The largest
``infinite crystal'', denoted by $K_{-\infty}^\infty$, arises when $c_i=\infty$
and $d_i=-\infty$ for all $i$. Among the variety of ``crystals'' produced by
the construction, $K_{-\infty}^\infty$ is distinguished by the property that
any monochromatic path in it is fully infinite, i.e., infinite in both forward
and backward directions (this object was introduced as the free combinatorial
A-type crystal by Berenstein and Kazhdan~\cite{BK-00}).

Equivalently: the principal lattice of $K_{-\infty}^\infty$ is formed by the
vertices $v[\bfa]$ for all $\bfa\in\Zset^n$. Also $K_{-\infty}^\infty$ can be
regarded as the ``universal'' RAN-crystal (with $n$ colors), due to the
following property:
  \begin{numitem}
any finite or ``infinite'' RAN-crystal is a (finite or infinite) principal
interval of $K_{-\infty}^\infty$, and vice versa.
   \label{eq:prop_inf}
  \end{numitem}
(An infinite principal interval of the form $\Int(v[\bfa],\infty)$ (resp.
$\Int(-\infty,v[\bfa])$) is the union of all paths beginning at $v[\bfa]$
(resp. ending at $v[\bfa]$).)

\section{Subcrystals with $n-1$ colors} \label{sec:subcryst}

In this section we apply the crossing model to study $(n-1)$-colored
subcrystals of an RAN-crystal $K=K(\bfc)=(V,E_1\cup\ldots\cup E_n)$.

For a subset $J\subset\{1,\ldots,n\}$ of colors, let $\Kscr(J)$ denote the
set of maximal connected subgraphs of $K$ whose edges have colors from
$J$, i.e., the components of the graph $(V,\cup(E_i\colon i\in J))$. When
the colors in $J$ go in succession, i.e., $J$ is an interval of
$(1,\ldots,n)$, each member $K'$ of $\Kscr(J)$ is a regular
$A_{|J|}$-crystal. (When $J$ has a gap, $K'$ becomes the Cartesian product
of several regular crystals. For example, for $J=\{1,3\}$, $K'$ is the
Cartesian product of two paths, with the color 1 and the color 3, or a
regular $A_1\times A_1$-crystal.)

We are interested in the case when $J$ is either $\{1,\ldots,n-1\}$ or
$\{2,\ldots,n\}$, denoting $\Kscr(J)$ by $\Kscr^{(-n)}$ in the former
case, and by $\Kscr^{(-1)}$ in the latter case. In other words,
$\Kscr^{(-n)}$ (resp. $\Kscr^{(-1)}$) is the set of $(n-1)$-colored
crystals arising when the edges with the color $n$ (resp. 1) are removed
from $K$.

Consider $K'\in\Kscr^{(-n)}$ and let $\Fscr(K')$ denote the set of
feasible functions corresponding to the vertices of $K'$. Since $K'$ is
connected, any $f\in\Fscr(K')$ can be obtained from any other
$f'\in\Fscr(K')$ by a series of forward and backward moves in levels
$1,\ldots,n-1$. So all functions in $\Fscr(K')$ have one and the same
tuple of values within the level $n$ of $G$. This level consists of nodes
$v^1_n(n),v^2_n(n-1),\ldots,v^n_n(1)$ (from right to left), and we denote
the $n$-tuple $(f(v^1_n(n)),\ldots,f(v^n_n(1)))$ by $\bfa(K')$, where
$f\in\Fscr(K')$. Thus, we have the following property: each subcrystal
$K'\in\Kscr^{(-n)}$ contains at most one principal vertex $v$ of $K$, in
which case $v=v[\bfa]$ for $\bfa=\bfa(K')$. Also the members of
$\Kscr^{(-n)}$ cover all principal vertices of $K$.

Similarly, for $K''\in\Kscr^{(-1)}$ and for the set $\Fscr(K'')$ of
feasible functions corresponding to the vertices of $K''$, the tuple
$\bfa(K''):= (f(v^1_1(1)),\ldots,f(v^n_1(1)))$ (where the nodes follow
from left to right in the level 1) is the same for all $f\in\Fscr(K'')$.
So each subcrystal $K''\in\Kscr^{(-1)}$ contains at most one principal
vertex of $K$ as well, and the members of $\Kscr^{(-1)}$ cover all
principal vertices of $K$.

We show a sharper property.
  \begin{prop}  \label{pr:subcryst}
Each subcrystal in $\Kscr^{(-n)}$ contains precisely one principal
vertex of $K(\bfc)$, and similarly for the subcrystals in
$\Kscr^{(-1)}$. In particular, $|\Kscr^{(-n)}|=|\Kscr^{(-1)}|=
(c_1+1)\times\ldots\times(c_n+1)$.
  \end{prop}

(This property need not hold when an $(n-1)$-element subset of colors is
different from $\{1,\ldots,n-1\}$ and $\{2,\ldots,n\}$.)

\smallskip
\begin{proof}
For a node $v$ of the supporting graph $G$, the maximal path beginning at
$v$ and going in the NE~direction is called the {\em NE-path} from $v$ and
denoted by $\PNE(v)$. Similarly, the maximal path beginning at $v$ and
going in the SE~direction is called the {\em SE-path} from $v$ and denoted
by $\PSE(v)$.

Let $K'\in\Kscr^{(-n)}$ and let $\bfa=(a_1,\ldots,a_n):=\bfa(K')$.
Consider an arbitrary function $f\in\Fscr(K')$. We show that the principal
function $f[\bfa]$ can be reached from $f$ by a series of forward moves,
followed by a series of backward moves, all in levels $\ne n$, whence the
desired inclusion $f[\bfa]\in\Fscr(K')$ will follow.

To show this, let $\Fscr_0$ be the set of functions $f'\in\Fscr(K')$ that
can be obtained from $f$ by (a series of) forward moves in levels $\ne n$
and such that $f'(v^k_k(1))=f(v^k_k(1))=: b_k$ for $k=1,\ldots,n$. Take a
maximal function $f_0$ in $\Fscr_0$. We assert that
  \begin{numitem}
the SW-edges of all nodes in $G$ (where such edges exist) are tight for
$f_0$.
   \label{eq:SWtight}
  \end{numitem}

Suppose this is not so for some node, and among such nodes choose a
node $v=v^k_i(j)$ with $i$ minimum. Acting as in
Section~\ref{sec:moves}, extend $G$ to the graph $\bar G$ and extend
$f_0$ to the corresponding function $\bar f_0$ on $V(\bar G)$ by
setting the upper bound $\bfb$ and the lower bound $\bfzero$ (then
$\bar f_0$ satisfies both the monotonicity condition and the switch
condition at each multinode and its values within each subgraph $G^k$
lie between 0 and $b_k$).

Consider an arbitrary node $v'=v_i^{k'}(j')$ with $1\le j'\le i$ in the
level $i$ and take the rhombus $\rho$ containing $v'$ as the right node;
let $z',u',w'$ be the left, upper and lower nodes of $\rho$, respectively.
Then $\partial \bar f_0(z',u')=0$ (this follows from $\bar f_0(z')=\bar
f_0(u')=c_{k'}$ when $j'=1$, and follows from the minimality of $i$ when
$j'>1$, in view of $(z',u')=\eSW(u')$). This implies $\eps(v')\ge
\delta(z')$ (where these numbers concern the bound $\bfb$);
cf.~\refeq{simple}. Moreover, this inequality is strict when $v'=v$ (since
$(w',v')=\eSW(v')$ and $\eSW(v)$ is not tight).

These observations imply $\tilde \eps_i(j)>0$, where $\tilde \eps$
concerns the bound $\bfb$. So the level $i$ contains an active multinode,
and therefore, $\bar f_0$ can be increased by a forward move in this
level. This move remains applicable when the bound changes to $\bfc$;
cf.~\refeq{eps-eps'}. Thus, $f_0$ is not maximal, and this contradiction
proves~\refeq{SWtight}.

From~\refeq{SWtight} it follows that for each $k$, all edges of the
NE-path from the bottommost node $\botvert^k$ in $G^k$ (going to the sink
$\rvert^k$) are $f_0$-tight. Hence $f(\rvert^k)=a_k$.

Now we apply (a series of) backward moves from $f_0$ in levels $\ne n$.
Let $\Fscr_1$ be the set of functions $f'\in\Fscr(K')$ that can be
obtained by such moves and satisfy $f'(\rvert^k)=a_k$ for $k=1,\ldots,n$.
Let $f_1$ be a minimal function in $\Fscr$. Arguing in a similar fashion,
one shows that
  \begin{numitem}
the NW-edges of all nodes in $G$ (where such edges exist) are tight for
$f_1$.
  \label{eq:NWtight}
  \end{numitem}
Now~\refeq{NWtight} implies that $f_1$ is constant within each $G^k$,
i.e., $f_1=f[\bfa]$, as required.

To show the assertion concerning $\Kscr^{(-1)}$, we can simply renumber
the colors, by regarding each color $i$ as $n-i+1$, and apply the model
for this numeration. Clearly the set of principal vertices preserves under
this renumbering, and now the result for $\Kscr^{(-1)}$ follows from that
for $\Kscr^{(-n)}$.
  \end{proof}

\noindent {\bf Remark 4.} Renumbering the colors as above causes a
``turn-over'' of the original model, so that level $i$ turns into level
$n-i+1$. (Note that the model does not maintain this transformation since the
switch condition~\refeq{feas}(iii) is imposed on SW- and SE-edges of nodes, but
not on NW- and NW-ones). A feasible function $f$ in the original model
corresponds to a feasible function $f'$ in the new model, so that $f$ and $f'$
determine the same vertex of the crystal. (In fact, the transformation
$f\mapsto f'$ is related to the Sch\"utzenberger involution in a crystal.) It
seems to be a nontrivial task to explicitly express $f'$ via $f$ (for $n=2$ an
explicit piece-wise linear relation is pointed out in~\cite{DKK-07}).

\medskip
We denote the member of $\Kscr^{(-n)}$ (of $\Kscr^{(-1)}$) containing a
given principal vertex $v[\bfa]$ by $K^{(-n)}[\bfa]$ (resp.
$K^{(-1)}[\bfa]$) and call it the {\em upper} (resp. {\em lower}) {\em
subcrystal of $K(\bfc)$ determined by $\bfa$}.

It turns out that one can explicitly express the parameter of
$K'=K^{(-n)}[\bfa]$.

To do this, note that the source and sink of $K'$ correspond to the
minimal function $f_{\rm min}(K')$ and the maximum function $f_{\rm
max}(K')$ in $\Fscr(K')$, respectively. One can see that in each $G^k$,
$f_{\rm min}(K')$ takes value 0 on all nodes, except for those on the path
$\PSE(\lvert^k)$ (from $\lvert^k$ to $\botvert^k$), where the values are
identically $a_k$. (The paths $\PNE(\cdot)$ and $\PSE(\cdot)$ are defined
in the proof of Proposition~\ref{pr:subcryst}.) In its turn, $f_{\rm
max}(K')$ takes value $c_k$ on all nodes, except for those on the path
$\PNE(\botvert^k)$ (from $\botvert^k$ to $\rvert^k$), where the value is
$a_k$. Symmetrically: the source and sink of a subcrystal
$K''=K^{(-1)}[\bfa]$ correspond to the minimal and maximum functions in
$\Fscr(K'')$, respectively, and in each $G^k$, the former takes value 0 on
all nodes, except for those on $\PNE(\lvert^k)$, where the value is $a_k$,
while the latter takes value $c_k$ on all nodes, except for those on
$\PSE(\topvert^k)$, where the value is $a_k$.
  \begin{prop}  \label{pr:parsubcryst}
The subcrystal $K^{(-n)}[\bfa]$ is isomorphic to the crystal
$K_{n-1}(\bfq)$ with colors $1,\ldots,n-1$, where $q_i=c_i-a_i+a_{i+1}$
for each $i$. In its turn, $K^{(-1)}[\bfa]$ is isomorphic to the crystal
$K_{n-1}(\bfq')$ with colors $2,\ldots,n$, where $q'_i=c_i-a_i+a_{i-1}$.
   \end{prop}
  \begin{proof}
Consider $f=f_{\rm min}(K^{(-n)}[\bfa])$ and $i\in\{1,\ldots,n-1\}$. From
the above description of $f$ it follows that for each node $v^k_i(j)$ with
$j\ge 1$ in the extended supporting graph $\bar G$, at least one of its
NW- and SW-edges is tight for $f$ (extended to $\bar G$), except possibly
for two nodes in the multinode $V_i(1)$: the first node $v=v_i^i(1)$, in
which $\partial f(\eSW(v))=\partial f(\eNW(v))=c_i-a_i$, and the second
node $v'=v_i^{i+1}(1)$, in which $\partial f(\eSW(v'))=a_{i+1}$ and
$\partial f(\eNW(v'))=c_{i+1}$. So, maintaining the monotonicity
condition~\refeq{feas}(i), one can increase the function (by the operator
$\phi_i$) only at $v$ or $v'$. More precisely, the active multinode in the
level $i$ is $V_i(1)$ (unless $q_i=0$) and the switch-node in it is either
$v$ or $v'$. If $a_{i+1}>0$, then $v$ cannot be the switch-node (since
$\eSW(v')$ is not tight). So the switch-node is $v'$, and the operator
$\phi_i$ acts $a_{i+1}$ times at $v'$, making the edge $\eSW(v')$ tight.
After that the switch-node becomes $v$ and $\phi_i$ acts $c_i-a_i$ times
at this node. This gives the desired parameter of $K^{(-n)}[\bfa]$.

(One can argue more formally. For each rhombus $\rho$ of $\bar G$ with the
left and right nodes in the level $i$, the value $\partial f(\rho)$
(defined in~\refeq{simple}) is zero, except possibly for two rhombi: the
rhombus $\rho$ whose right node is $v$, where $\partial f(\rho)=c_i-a_i$,
and the rhombus $\rho$ whose right node is $v'$, where $\partial
f(\rho)=a_{i+1}$. This implies that the total residual upper slack
$\tilde\eps_i(1)+\ldots+\tilde\eps_i(i)$ for $f$ in the level $i$ is just
$q_i=c_i-a_i+a_{i+1}$.)

The assertion concerning the lower subcrystal $K^{(-1)}[\bfa])$ follows by
symmetry (when each color $i$ is renumbered as $n-i+1$).
  \end{proof}

\noindent {\bf Remark 5.} This proposition implies that  all possible
parameters $\bfq$ of the upper subcrystals of $K(\bfc)$ give the set of integer
points of some polytope in $\Rset^{n-1}$. Note also that for corresponding
tuples $\bfq$ and $\bfa$, the numbers $a_1,\ldots,a_{n-1}$ are determined by
$\bfq$ and $a_n$, namely: $a_i=c[i:n-1]-q[i:n-1]+a_n$ for $i<n$. This enables
us to compute the quantity $\eta(q)$ of crystals in $\Kscr^{(-n)}$ having a
prescribed parameter $\bfq$: this is as large as the set of numbers
$a_n\in\Zset$ that together with $\bfq$ determine $\bfa$ satisfying  $0\le
a_i\le c_i$ for all $i=1,\ldots,n$. (One can express $\eta(q)$ as the
difference between $\min\{c_n,\, q[i:n-1]-c[i+1:n-1] \colon i=1,\ldots,n-1\}$
and $\max\{0,\,q[i:n-1]-c[i:n-1]\colon i=1,\ldots,n-1\}$. In particular, if
$c_i=0$ takes place for some $i$, then all upper subcrystals of $K(\bfc)$ are
different.) This gives a {\em branching rule} for decomposing an irreducible
$sl_{n+1}$-module into the sum of irreducible $sl_n$-modules. In the above
expression, the branching rule looks simpler than the rule indicated
in~\cite[Corollary~2.11]{BZ-01}.

 \medskip
Next, we are able to indicate where a principal vertex $v[\bfa]$ of $K(\bfc)$
is located in the subcrystals $K^{(-n)}[\bfa]$ and $K^{(-1)}[\bfa]$.
  \begin{prop} \label{pr:location}
Let $\Pi'$ be the principal lattice of the upper subcrystal
$K'=K^{(-n)}[\bfa]$. Then the principal vertex $v=v[\bfa]$ of $K=K(\bfc)$
is contained in $\Pi'$ and the $(n-1)$-tuple $\bfa'$ of its coordinates in
$\Pi'$ satisfies $a'_i=a_{i+1}$ for $i=1,\ldots,n-1$. Symmetrically, $v$
is contained in the principal lattice $\Pi''$ of the lower subcrystal
$K^{(-1)}[\bfa]$ and the $(n-1)$-tuple $\bfa''$ of its coordinates in
$\Pi''$ satisfies $a''_i=a_{i-1}$, $i=2,\ldots,n$.
  \end{prop}
  \begin{proof}
We have to show that the principal vertex $v'$ of $K'$ with the
coordinates $\bfa'$ in $\Pi'$ coincides with $v$. By explanations in
Subsection~\ref{sec:princ}.2 (applied to $K'$ in place of $K$), the vertex
$v'$ can be obtained from the source of $K'$ by applying the sequence of
forward moves (in $\Mscr_{n-1}$) corresponding to the combined string
  $$
  (S_{n-1,n-1})^{a_n} (S_{n-1,n-2})^{a_{n-1}}\cdots (S_{n-1,2})^{a_3}
  (S_{n-1,1})^{a_2}.
  $$

For $k=1,\ldots,n$, partition $V(G^k_n)$ into two subsets $L_k,R_k$, where
$L_k$ is the set of nodes of the path $\PSE(\lvert^k)$ and $R_k$ is the
rest. Note that $R_1=\emptyset$ and that for $k>1$, the subgraph of
$G^k_n$ induced by $R_k$ is isomorphic to $G^{k-1}_{n-1}$. Also: (a) the
minimal feasible function $f_{\rm min}(K')$ for the subcrystal $K'$ (in
$\Mscr_n$) takes the constant value $a_k$ on $L_k$ and 0 on $R_k$, for
each $k$; and (b) the principal function $f[\bfa]$ takes the value $a_k$
on each $L_k\cup R_k$.

Suppose that $\bfb\in\Zset_+^n$ is such that $\bfb\le \bfa$, and that $f$
is the feasible function taking the constant values $a_k$ and $b_k$ within
$L_k$ and $R_k$, respectively, for each $k$. To obtain the desired result,
it suffices to show the following:
  \begin{numitem}
let $b_k<a_k$ for some $k>1$, and let $f'$ be the feasible function (in
$\Mscr_n$) taking the constant value $b_k+1$ within $R_k$ and coinciding
with $f$ on the rest; then $f'$ is obtained from $f$ by applying the
sequence of moves (in $\Mscr_n$) corresponding to the fundamental string
$S_{n-1,k-1}$.
  \label{eq:bk}
  \end{numitem}

According to~\refeq{fstring}, $S_{n-1,k-1}=w_{n-1,k-1,n-k+1}\cdots
w_{n-1,k-1,2} w_{n-1,k-1,1}$, and for each $i$, the substring
$w_{n-1,k-1,i}=:w'_i$ is $(i)(i+1)\cdots (i+k-2)$. Observe that each
$w'_i$ corresponds to the (maximal) NE-path $P_i$ in $G^{k}_n$ beginning
at the node $v^k_{i+k-2}(i)$ (which is the $i$-th node on the SW-side of
the rectangular indiced by $R_k$, viz. on the path $\PSE(v^k_{k-1}(1))$).
One can check that the action corresponding to $S_{n-1,k-1}$ changes $f$
only within the base subgraph $G^k_n$ and the action corresponding to
$w'_i$ consecutively increases the current function along the path $P_i$.
This results in the function $f'$ as required in~\refeq{bk}. A
verification in details is left to the reader.

The assertion concerning $K^{(-1)}[\bfa]$ follows by symmetry.
  \end{proof}

Note that for $K'=K^{(-n)}[\bfa]$ and $\bfa'$ as in
Proposition~\ref{pr:location}, if we apply the first part of
Proposition~\ref{pr:parsubcryst} to $K'$ and the second part to the lower
subcrystal $\tilde K$ of $K'$ determined by $\bfa'$, then we obtain that
the $(n-2)$-colored crystal $\tilde K$ with colors $2,\ldots,n-1$ that
contains the principal vertex $v[\bfa]$ of $K(\bfc)$ has the parameter
$\bfc''$ such that $c''_i=(c_i-a_i+a_{i+1})-a'_i+a'_{i-1}=
c_i-a_i+a_{i+1}-a_{i+1}+a_i=c_i$ for $i=2,\ldots,n-1$. (This $\tilde K$ is
the component of $K^{(-n)}[\bfa]\cap K^{(-1)}[\bfa]$ that contains
$v[\bfa]$.) This leads to a rather surprising property:
  \begin{numitem}
for any $1\le r\le\lceil n/2\rceil$, all $(n-2r+2)$-colored subcrystals of
$K(\bfc)$ with the colors $r,r+1,\ldots,n-r+1$ that meet the principal
lattice $\Pi(\bfc)$ have the same parameter, namely,
$(c_r,\ldots,c_{n-r+1})$, and therefore, they are isomorphic.
  \label{eq:par_middle}
  \end{numitem}

Finally, using the crossing model, one can compute the lengths of maximal
monochromatic paths in $K^{(-n)}[\bfa]$ (or in $K^{(-1)}[\bfa]$) that go
through the principal vertex $v[\bfa]$ of $K$ (one can say that the length
concerning a color $i$ expresses the ''$i$-width" of the subcrystal at
this vertex).

 \medskip
In conclusion of this paper we can add that the crossing model can be used to
reveal more structural properties of RAN-crystals. A nontrivial problem on this
way is to characterize the intersection of the upper subcrystal
$K^{(-n)}[\bfa]$ and the lower subcrystal $K^{(-1)}[\bfb]$ of $K(\bfc)$ for any
$\bfa,\bfb\in\Zset_+^n$ (this intersection may be empty or consist of one or
more subcrystals with colors $2,\ldots,n-1$). This problem is solved in the
forthcoming paper~\cite{DKK-08b}, giving rise to an efficient recursive
algorithm of assembling the RAN-crystal $K(\bfc)$ for a given parameter $\bfc$.
Also using the model, we explain there that a regular $B_n$-crystal
($C_n$-crystal) with parameter $\bfc=(c_1,\ldots,c_n)$ can be extracted from
the ``symmetric part'' of the regular $A_{2n-1}$-crystal with parameter
$(c_1\ldots,c_n,\ldots,c_1)$ (resp. the regular $A_{2n}$-crystal with parameter
$(c_1,\ldots,c_n,c_n,\ldots,c_1)$).

 \smallskip
{\bf Acknowledgements.} We thank the anonymous referees for remarks and
useful suggestions.

\small
\bibliographystyle{plain}

\end{document}